\documentclass[11pt]{article}
\usepackage{a4wide}
\usepackage{latexsym}
\usepackage{amsmath,amsthm,amsfonts,amssymb}
\usepackage{mathrsfs}
\usepackage{wasysym}
\usepackage{amsfonts}
\usepackage{amssymb}
\usepackage{epsfig}
\usepackage[all]{xy}
\newtheorem{defi}{\textbf{Definition}}[section]
\newtheorem{theo}[defi]{\textbf{Theorem}}
\newtheorem{lemma}[defi]{\textbf{Lemma}}
\newtheorem{prop}[defi]{\textbf{Proposition}}
\newtheorem{defprop}[defi]{\textbf{Definition/Proposition}}
\newtheorem{coro}[defi]{\textbf{Corollary}}
\newcommand{\Ang}{{\rm Ang}}
\newcommand{\length}{{\rm length}}
\newcommand{\Cone}{{\rm Cone}}

\newcommand{\calL}{\mathcal{L}}
\newcommand{\calC}{\mathcal{C}}
\newcommand{\tilCaG}{\widehat{Cay}\Gamma}
\newcommand{\frap}{\mathfrak{p}}
\newcommand{\G}{\Gamma}
\newcommand{\tilG}{\tilde{\Gamma}}

\newcommand{\oup}{(\Omega \cup P)}
\newcommand{\tte}{{\mathtt{eq}}}
\newcommand{\ttin}{{\mathtt{ineq}}}
\newcommand{\ttco}{{\mathtt{cons}}}
\newcommand{\calR}{\mathcal{R}}
\newcommand{\calW}{\mathcal{W}}
\newcommand{\calV}{\mathcal{V}}
\newcommand{\calD}{\mathcal{D}}
\newcommand{\calX}{\mathcal{X}}
\newcommand{\calS}{\mathcal{S}}
\newcommand{\calT}{\mathcal{T}}

\newcommand{\absx}{|_{\mathcal{X}}}
\newcommand{\dron}{\partial}
\newcommand{\tilHi}{\tilde{H}_i}

\newcommand{\tto}{\twoheadrightarrow}
\newcommand{\Sec}{{\mathrm{Sec}}}

\title{Existential questions in (relatively) hyperbolic groups {\it and} Finding relative hyperbolic structures}
\author{Fran\c{c}ois Dahmani}

\date{}
\begin{document}

\maketitle
{\footnotesize \begin{center}
{\bf Abstract. } This arXived paper has two independant parts, that are improved and corrected versions of different parts of a single paper, once named ``On equations in relatively hyperbolic groups''.

The first part is entitled ``Existential questions in (relatively) hyperbolic groups'' (to appear in {\it Israel J. Math.}).  We study there the existential theory of  torsion free 
hyperbolic and relatively hyperbolic groups, 
in particular those with virtually abelian
parabolic subgroups. We show that the satisfiability of systems of equations
and inequations is decidable in these groups.

In the second part, called ``Finding relative hyperbolic structures'' (to appear in {\it Bull. London Math. Soc.}), we provide a  general algorithm that recognizes the  class of groups that are  hyperbolic relative to abelian subgroups.

\end{center}

}

\bigskip
\part{Existential questions in (relatively) hyperbolic groups}

The existential theory of a group $G$ is the set of all sentences of first
order logic, in the language of groups, that contain only existential
quantifiers $\exists$, and that evaluate to true in $G$.  
Each such sentence is
equivalent to the existence of a solution of the disjunction of finitely many
finite systems of equations (of the form $x_1 x_2 \dots x_k =1$) and
inequations (of the form $y_1 y_2 \dots y_n  \neq 1$), with parameters in the
group.  
Thus deciding whether a sentence in the existential theory of a group
is true amounts to deciding whether a system of equations and inequations
(with parameters) admits a solution in the group.  
If an algorithm can
perform this task, we say that the group under consideration has decidable existential theory with parameters. 
Sometimes, one is more interested in solving finite systems of equations, with parameters, without considering inequations, and one then speaks of  the Diophantine theory.  
Also commonly used is the universal theory of a group, that is the set of sentences of first order logic
 that use only universal quantifier $\forall$. 
The latter is decidable if and only if the existential theory is, 
 since one goes from one to the other by negations.

To have decidable Diophantine theory, or even better, existential theory, is a
very strong property.  
For example, the Conjugacy Problem (deciding whether
two given elements of a group are conjugate to each other) is only a single
instance of the Diophantine theory. 
The study of existential and Diophantine
theory  of a group $G$ is also closely related to the study of morphisms from finitely
presented groups to  $G$, since each  morphism $H\to G$ is a
solution to a system of equations in $G$: the unknowns are the images of the
generators of $H$ and the equations in $G$ are the defining relations of $H$.

  There are rather few classes of groups that are known to have decidable
  existential theory. 
Clearly this is true for each finite group, and abelian
  group. 
For free groups, this is a celebrated  result of G.~Makanin
  \cite{Maka}. 
Recently, V.~Diekert and A.~Muscholl proved the decidability
  for right-angle Artin groups \cite{DiMu}, and more generally, V.~Diekert and
  M.~Lohrey \cite{DiLo} proved the stability of the decidability of
  existential theory under graph products (in particular free products, and
  direct products). 
For torsion free hyperbolic groups, the problem has been
  open for some time, since their Diophantine theory was proved decidable by
  E.~Rips and Z.~Sela \cite{RS}.   
The question indeed appears in Sela's list
  of problems \cite{Slist}. 
In fact, he gave  himself an answer  in the very last pages of his
  seventh paper \cite{SVII} of his series about elementary theory of free and
  hyperbolic groups, \cite{SVII}, using deep results about strict
  Makanin-Razborov resolutions. 
We propose an alternative proof here, possibly simpler, 
  although still relying on Rips and Sela's method for the
  Diophantine theory.

Solving equations in groups is an important and difficult problem, and has
many applications. 
But the problem of ``solving'' can be interpreted in
several ways, for instance in finding the algebraic structure of the set of
solutions, or in algorithmically deciding the existence of solutions. As we
said, we are interested here in the second (weaker) interpretation, but we
believe (and this is illustrated in further work \cite{DG}) that our methods
allow to tackle the first, to some extent.

Thus, the general problem addressed in this paper is to find, for certain groups, an algorithm that, given a finite system of equations and inequations with parameters in the group, indicates whether it admits a solution. 
The groups we are interested in are  hyperbolic, and relatively hyperbolic groups. 
 Both classes  were introduced by M. Gromov \cite{Grom} and the theory of relatively hyperbolic groups was developed by B. Farb \cite{F}, B. Bowditch \cite{Brel}, and recently by many authors. 
These groups are coarse analogues of the  fundamental groups of finite volume hyperbolic manifolds. The latter are hyperbolic relative to the collection of their cusps subgroups. In this point of view, we say that a group is hyperbolic relative to a collection of subgroups (called its parabolic subgroups) if it acts on a proper Gromov-hyperbolic space, preserving a collection of disjoint horoballs, cocompactly on the complement of these horoballs, and such that the stabilizers of the horoballs, are exactly the parabolic subgroups. 

In principle, relatively hyperbolic groups can have parabolic subgroups of arbitrary algebraic nature. 
However, in my opinion, the most natural examples of relatively hyperbolic groups have abelian or virtually abelian parabolic subgroups. 
This is the case for fundamental groups of finite volume (or even geometrically finite) hyperbolic (or even sufficiently pinched negative curvature) manifolds, for limit groups, or fully residually free groups, as I observed in previous work \cite{Dgt}, as well as for the larger class of groups acting freely on $\mathbb{R}^n$-trees (for the same reason actually), and for co-compact isometry groups of $CAT(0)$ spaces with isolated flats (as proved by B.~Kleiner and C.~Hruska \cite{KH}).  
Of course it is even the case for hyperbolic groups.

We prove:

\begin{theo}\label{theo;intro2}
 Let $\Gamma$ be a torsion-free group, hyperbolic relative to a family of 
 subgroups that have decidable existential theories with parameters. 
 Then the existential theory with parameters of  $\Gamma$ is decidable. 

 In particular, the  existential theory with parameters of any torsion free hyperbolic group is decidable.
\end{theo}

 For the latter assertion, we also refer to \cite[Theorem 7.12]{SVII}.

  Theorem \ref{theo;intro2} immediately implies that torsion-free groups that are hyperbolic relative to free abelian subgroups 
  have  decidable existential theory. 
A special case is the study of O.~Kharlampovich, A.~Miasnikov's 
 of satisfiability of equations in fully residually free (or limit) groups \cite{KM}, which we thus recover in a possibly
 simpler way. 
 However their result is stronger in the direction of systems of equations, since it gives the structure of solutions in
 the form of Makanin-Razborov diagrams.

  A folklore result states that  virtually abelian groups have decidable
  existential theory.  
We did not know any reference for that, so we provide a
  proof here. It follows that:

\begin{coro}\label{coro;intro3}
   The existential theory with parameters of any torsion free relatively hyperbolic group with virtually abelian 
   parabolic subgroups, is decidable.
\end{coro}

This answers, in particular, a question of Sela in his problem list \cite{Slist} about the Diophantine theory of $CAT(0)$ groups with isolated flats.

The ``next'' most natural class of relatively hyperbolic groups is probably
that of groups with virtually nilpotent parabolic subgroups, since it contains all pinched negatively curved finite volume manifolds. Unfortunately, there seems to be little hope in this direction. 
In fact, Roman'kov  has produced a 4-step nilpotent group with undecidable existential theory \cite{Ro}. 
His proof  uses  Yu.~Matijasevic's solution to the  tenth Hilbert problem on the un-decidability of the existential theory of the ring of integers.

In this work, only the torsion free case has been considered. At the end of the
first part, we give a brief account of the difficulties occurring in presence
of torsion. 
  The general approach follows that of Rips and Sela: to use canonical
  representatives, in order to lift equations in a hyperbolic group into
  equations in a free group. 
For equations in a relatively hyperbolic group, we will use an appropriate free product instead of a free group. 
The main new idea, that allows to solve
  inequations, is to consider rational constraints, that are
  requirements that some unknowns lie in a language recognized by a finite
  state automaton, in some way. 
We in fact simulate the inequations in a
  hyperbolic group by rational constraints in the free group into which the
  system of equations is lifted. 
In general such a tactic has little chance to
  work, since the kernel of a quotient of a free group is never rational
  unless it is of finite index. 
However, our study can be done by first
  restricting the search of lifts to a specific  subset of the free group, or
  the free product, that will be proven rational. 
This is the technical part
  of the paper.

  All the procedure  is first illustrated in the easier case of hyperbolic
  groups in the first section. In this context many technical difficulties are
  avoided, and the proof is much shorter.

It should be noted that the algorithms we provide are explicit and can be (theoretically) implemented provided the following technical items are known or given: a presentation of the ambient group, generators and defining relations for each parabolic subgroup up to conjugacy,  
a hyperbolicity constant of the associated coned-off graph, the list of simple loops of any given length up to some explicit length $\kappa$, up to translation  in this graph, 
 and finally  machines that solve the existential theories of the parabolic subgroups. Thus, the computability of different tools involved in this work is to be understood with this preliminary knowledge given, or found. 
It might be interesting to note that, if the parabolic subgroups are abelian,  we reduce this in part 2 (Coro. 7.4, or 2.4 in the published paper), to the single preliminary knowledge of a finite presentation.

I wish to thank  M.~Sapir,  O.~Kharlampovich and T.~Delzant, for their encouragements, and    I.~Bumagin, V. Guirardel, and D.~Groves for interesting related discussions. My acknowledgments also go to the referee for some useful remarks.

   \section{Existential Ideas}

      The aim of this first section is to introduce notations and formalism of equations inequations and constraints, 
      and to expose the ideas of the paper in the special, but motivating, case of torsion free hyperbolic groups.

      As said in the introduction,  the \emph{existential theory} of a group is the set of all first order
      logic sentences that use only the existential quantifier $\exists$, 
      and that are  satisfied by the group. Any such sentence is equivalent to
     a disjunction of finitely many pairs 
      of systems of equations and of inequations. Thus
       deciding whether an existential sentence is true in a group is
      equivalent to deciding whether a system of equations and inequations has
      a solution.

      It is worth noting that we allow the presence of parameters (or
      constants) in our equations. 
      This is standard for systems of equations alone 
      (indeed any system of equations without parameter has a trivial
      solution), but not always for systems of equations and inequations.

      \subsection{Formalities about equations.}
         Let us start with a group $\G$. 
         Let $\Omega$ be an abstract finite set,  
	that we will call the set of
         unknowns, 
         and let $P$ be a finite subset of $\G$  
         that we will call the set of parameters.

         \paragraph{Systems.}

            A \emph{system of equations} in $\G$, with unknowns 
            in $\Omega$, and parameters in $P$, is a family of expressions 
            of the form $\tte=1$, where $\tte$ is an element of $\oup^*$, that
            is a word with letters in $\oup$. 
            Let us choose notations for the future, and denote a  finite system of equations by 
            $\mathscr{E}=\{\tte_1=1, \tte_2=1, \dots, \tte_n=1  \}$,  where each $\tte_i$ is in $\oup^*$. 
	   Let us remark that one usually considers equations where inverses of the unknowns appear. 
	   This can easily be reduced to the setting above by adding an equation $xy=1$ and using $y$ 
		instead of $x^{-1}$ everywhere else.

            A system of \emph{inequations} in $\G$, is a family of  expressions 
            of the form $\ttin \neq 1$, where $\ttin$ is an element of $\oup^*$. 
            Similarly let us denote a  finite system of inequations by 
            $\mathscr{I}=\{\ttin_1\neq 1, \ttin_2\neq 1, \dots, \ttin_m\neq 1
            \}$,  where each $\ttin_i$ is in $\oup^*$.

            There is a refinement to the concept of inequations, that is the concept of constraint. 
            A system of \emph{constraints} in $\G$, is a family of  expressions 
            of the form $\ttco \in L$, where $\ttco$ is an element of  $\oup^*$, and $L$ is a subset of $\G$. 
            For the moment, we do not wish to add restriction on what kind of
            subset it is, this topic will be developed later on.
            Let us denote a finite system of constraints by $\mathscr{C}=\{\ttco_1\in L_1, \dots, \ttco_p \in L_p \}$.
            
         \paragraph{Solutions, and satisfiability.}

            Any application $s: \Omega \to \G$,   
	induces a morphism of monoids 
            $s: \oup^* \to \G$, simply by deciding that $s$ is  the identity
            on $P$.

            A \emph{solution} to the system $(\mathscr{E},\mathscr{I}, \mathscr{C} )$ 
            is an application $s:\Omega \to \G$ such that the induced morphism 
            satisfies $ \forall i\leq n, \, s(\tte_i) = 1$ in $\G$, 
            $\forall j\leq m, \, s(\ttin_j) \neq 1$ in $\G$, and  
            $s(\ttco_k) \in L_k, \forall k\leq p$.

            We say that a system of equations, 
            with possibly inequations and constraints, is \emph{satisfiable}, 
            if there exists at least one solution.

         \paragraph{Triangulation.}

            With an easy triangulation procedure,  
            given a system $(\mathscr{E},\mathscr{I}, \mathscr{C} )$ 
            with unknowns $\Omega$,   one can construct another system 
            $(\mathscr{E}',\mathscr{I}', \mathscr{C}')$, 
		satisfiable if and only if  
            $(\mathscr{E},\mathscr{I}, \mathscr{C} )$ is satisfiable,  
            such that each constraint $\ttco' \in \mathscr{C}'$ and each
            inequation  $\ttin' \in \mathscr{I}'$ is a word of only 
		one letter in $\Omega'$, 
            and such that each equation $\tte' \in \mathscr{E}'$ has length $3$
            as a word in $(\Omega' \cup P)^*$ (equations of length $2$ can be made of length $3$ by multiplying with the neutral element of the group, seen as a parameter in $P$, or alternatively, one can discard them using the assumption of absence of element of order $2$ in the group).

		Of course, this triangulation may increase 
            the number of unknowns, and of equations. 

            We will always assume that this reduction has been done.

      \subsection{Hyperbolic groups.}

         We now turn to a more informal discussion of the case of hyperbolic groups.
         Let $\G$ be a torsion free hyperbolic group. The problem of deciding the 
         satisfiability of equations (with parameters) in $\G$ was solved in
         the positive way by E.~Rips and Z.~Sela in \cite{RS}. 
         
         Their strategy is as follows. Given a  system of equations
         $\mathscr{E}$, they transform it  into a system of equations whose 
         unknowns are paths in the Cayley graph. This is interpreted as
         equations in the free group over the generators of $\G$.

         In their construction, a solution of the new equations of paths always  produces  a solution of $\mathscr{E}$. 
         The converse is by no way easy, and is the object of the theorem of
         canonical representatives which we state now, because we will need the constants $\kappa,\lambda,\mu$ for the construction.

           \begin{theo}(\cite{RS}, Canonical representatives) \label{theo;repcanhyp}

             Let $\Gamma$ be a torsion free hyperbolic group, and $S$ a finite
             symmetric generating set. 
             There exists computable 
             constants $\kappa, \lambda, \mu>0$, such that the following holds.

             Let $\mathcal{A}$ be an arbitrary finite subset of $\G$, $n$ be an integer, and $\tte_1, \dots,
             \tte_n$  be $n$
             words of three letters in  
             $\mathcal{A}$, such that, for all $i$, $\tte_i  \underset{\G}{=}   1$ in
             $\Gamma$. Let us note $\tte_i=a(i)_1 a(i)_2 a(i)_3$ for all $i$.

             Then, for all $a\in \mathcal{A}$, there exists $\tilde{a}\in F_S$, labeling a 
             $(\lambda,\mu)$-quasi-geodesic in $Cay \G$ from $1$ to $a$, such that 
             \[ \forall i, \;  \exists l_1,l_2,l_3,   \;  \exists r_1,r_2,r_3, \; \exists c_1,c_2,c_3 \]
             satisfying   
             \[  \forall j\leq 3, \quad |c_j| \leq \kappa n, \quad
             \widetilde{a(i)_j} = l_j c_j r_j \quad r_j = (l_{j+1
              \, [3]})^{-1}. \]

           \end{theo}

           We call the elements $\tilde{a}$ the canonical representatives of
           the elements of  $\mathcal{A}$.  In \cite{RS}, the authors introduce a number of candidates, depending on a parameter, to be canonical representatives, and prove that for a certain parameter, they satisfy the relations stated. But there seems to be no way to find the parameter {\it a priori}. We are not interested here in that computability, but only on the existence of such representatives.

         \paragraph{Path equations.}

           We explain now how to transform a system of equations in $\G$
           into a family of systems of 
           equations of paths in $Cay \G$ 
           (for a certain finite symmetric generating set $S$ of $\G$).

           First, let us denote by $F_S$ the free group over $S$. 
           Each element of $F_S$ has a 
           normal form which is a word over the alphabet $S$, 
           and thus labels a unique path in $Cay \G$ starting at $1$. 
           Conversely, for any path in $Cay \G$, the word read on the
           consecutive edges of the path defines a unique element of $F_S$.    
           If $w\in F_S$, let us denote by $\bar{w}$ its image in $\Gamma$. 

           Let us write each equation in $\mathscr{E}$ as $\tte_i = z(i)_1
           z(i)_2 z(i)_3$, where $z(i)_j$ is an element of $\oup$.

	  Let $\calC_0$ be the set of all possible $n$-tuples of triples $ ( (c(i)_1,c(i)_2,c(i)_3) )_{i\leq n}$ in which $c(i)_j$ is an element of  $F_S$ of length at most $\kappa n$ satisfying    $\overline{c(i)_1}\overline{c(i)_2}\overline{c(i)_3}   \underset{\G}{=}  1$. Since $\kappa$ is an explicit constant depending on $\delta$ and the group presentation,  this is a finite computable set.

	For all parameter $p\in P$, Let $\calC_p$ be the set of all 	elements of $F_S$ labeling a $(\lambda,\mu)$-quasi-geodesic from $1$ to $p \in \Gamma$. This is again a computable set.

	Now for all choice of $c \in \calC_0$ and of $rep(p) \in \calC_p$, we write the system  $\mathscr{E}'(c,rep)$ in $F_S$:

           \[
           \begin{array}{lclcl}
             \widetilde{z(i)_1} = {l(i)_1}{c(i)_1}{r(i)_1}  &                 & r(i)_1 = l(i)_2^{-1}   \\
             \widetilde{z(i)_2} = {l(i)_2}{c(i)_2}{r(i)_1}  &  \hbox{       } & r(i)_2 = l(i)_3^{-1}   \\
             \widetilde{z(i)_3} = {l(i)_3}{c(i)_3}{r(i)_1}  &                 & r(i)_3 = l(i)_1^{-1}
           \end{array}
           \]

           The $c(i)_j$ are parameters (given by the choice $c\in\calC$) that satisfy 
           $\overline{c(i)_1}\overline{c(i)_2}\overline{c(i)_3}   \underset{\G}{=}  1$. 
           The set of unknowns is $\tilde{\Omega} \cup \{ l(i)_j, r(i)_j,
           i\leq n, j\leq 3\}$, where  $\tilde{\Omega} \simeq \Omega$. 
           If $z(i)_j \in \Omega$ then $\widetilde{z(i)_j}$ has the
           corresponding value in $\tilde{\Omega}$, 
           and if $z(i)_j \in P$, 
           then  $\widetilde{z(i)_j}= rep(z(i)_j)$ is our chosen parameter in $F_S$.

		\begin{prop}\label{prop;equiv_syst_hyp}
			The system $\mathscr{E}$ has a solution in $\G$ if and only if, for some choice of $c$ and $rep(p), p\in P$ as above,  
			the system  $\mathscr{E}'(c,rep)$ has a solution.  
		\end{prop}

           It is clear that any solution $s$ of any $\mathscr{E}'(c,rep)$ provides a solution for $\mathscr{E}$,   by taking the images in $\G$ of  the values of $s(\tilde{\Omega})$, 
           that is by  considering 
           $\bar{s}$ on $\Omega$: $\bar{s}(\omega) = \overline{s(\tilde{\omega})}$.

           The converse is the Theorem \ref{theo;repcanhyp}. $\square$

         \paragraph{Satisfiability of equations.}

           This is already sufficient to decide the satisfiability of $\mathscr{E}$ in $\Gamma$. 
           Indeed, this is now equivalent to the satisfiability of one of the systems $\mathscr{E}'(c,rep)$ in $F_S$. 
           Recall that these systems were effectively computed from  $\mathscr{E}$. 
           We can  use known solutions to 
           the satisfiability problem for equations in free groups, 
           such as those proposed by Makanin \cite{Maka},  
           or Diekert-Guti\'errez-Hagenah \cite{DiGuHa}.
           
           We have finished our review of Rips and Sela's result on equations
           in 
           torsion free hyperbolic groups. 
           We now go further, and explain how to deal with inequations.

         \paragraph{Geometric solutions.}

             Theorem \ref{theo;repcanhyp} says more than already
            used. 
            It allows to lift any solution of $\mathscr{E}$ into a solution of
            one of the $\mathscr{E}'(c,rep)$, which has the particularity that it consists of
            elements that label quasi-geodesic paths in $Cay \G$.  Many
            solutions to  $\mathscr{E}'(c,rep)$ probably do not have this property,
            and thus are of little relevance for our study of
            $\mathscr{E}$. 

            In particular, if one wants to consider some
            inequations in $\Gamma$ in addition of  $\mathscr{E}$, then it
            will be hard to tell which solutions of  $\mathscr{E}'$ actually
            map on solutions of $\mathscr{E}$ that satisfy the inequation. But
            if one considers only solution with a geometric relevance, this
            control is made easier.

            Recall that a $L$-local $(\lambda,\mu)$-quasi-geodesic is an arc-length
            parameterized  path such that any subsegment of length at most $L$
            is a $(\lambda,\mu)$-quasi-geodesic.
            We use the classical result of hyperbolic geometry, that for all
            $(\lambda,\mu)$, there exists $L,\lambda',\mu'$ so that any 
            $L$-local $(\lambda,\mu)$-quasi-geodesic is a global  
            $(\lambda',\mu')$-quasi-geodesic (see \cite[Chap. 3, Thm 1.4]{CDP}, including 
	   explicit values of these constants in term of $\delta$). 

            We then call \emph{geometric} an element of $F_S$
            that labels a $L$-local $(\lambda,\mu)$-quasi-geodesic in $Cay
            \G$. Let us denote by $\calL$ the set of geometric elements in 
            $F_S$. 

            Every canonical representative is geometric, and
            therefore, the equivalence of 
		Proposition \ref{prop;equiv_syst_hyp} can be refined: 
	    the solutions of  $\mathscr{E}$
            are exactly the images of 
            the geometric solutions of $\mathscr{E}'(c,rep)$ for $c\in \calC_0$ and $rep(p)\in \calC_p\,  \forall p \in P$.

         \paragraph{Inequations in $\G$ as constraints in geometric solutions.}

         It is then easy to translate what it means to satisfy inequations 
         in $\G$. Indeed, the geometric elements of $\G$ that define the
         trivial element of $\G$ form a finite computable set, since they are
         all of length at most $\mu'$. Let $\calL_0$ be this set.

         Let  $\mathscr{I}=\{\ttin_1\neq 1, \ttin_2\neq 1, \dots, \ttin_m\neq 1
            \}$, be a system of inequations where each $\ttin_i$ is in
            $\Omega$.

            Then the system $(\mathscr{E}, \mathscr{I})$ has a solution if and
            only if one of the $\mathscr{E}'(c,rep)$ has a solution $s$ consisting of geometric
            elements, and such that  each $s(\widetilde{\ttin_i})$, 
            for $i\leq m$, 
            in outside $\calL_0$. 
         
            In other words, we have the proposition:

            \begin{prop}\label{prop;equiv}
              Let $(\mathscr{E}, \mathscr{I})$ be a system of equations and
              inequations in $\G$, triangulated as above. Let $\mathscr{E}'(c,rep)$
              be the finite family of systems of equations in $F_S$ constructed above. 

              Let us define the  constraints in $F_S$,  
            $\mathscr{C}'=\{\tilde{\omega} \in \calL, \forall \tilde{\omega} 
            \in \tilde{\Omega}
            \} \cup \{\widetilde{\ttin_i} \in \calL \setminus \calL_0, \forall
            i\leq m\}$. 

            Then the system $(\mathscr{E}, \mathscr{I})$  has a
            solution in $\Gamma$ if and only if one of the systems $ (\mathscr{E}'(c,rep),\mathscr{C}')$
            has a solution in $F_S$.

            \end{prop}

         \paragraph{Satisfiability.}
            The interest of this transformation is that the satisfiability of 
            the system  $(\mathscr{E}'(c,rep),\mathscr{C}')$ in $F_S$ is decidable.

 The decidability of existential theory for free groups is originally
      due to G.~Makanin \cite{Maka}, who previously solved the case of 
      ``word'' equations,
      in free monoids (without involution), in the 80's. 
      K.~Schulz introduced treatment of constraints in his
      algorithm for word equations, in the 90's. 
      An alternative  solution for word equations, with rational constraints, 
	was found by  
      Plandowski. This
      has  advantages about complexity, and could be also used for
      equations with constraints in free groups.

            A rational language in a free group is a subset whose set of
            normal forms is  a word language
            recognized by a finite state automaton. A rational constraint is
            a constraint defined by a rational language.

      \begin{theo}(Diekert, Guti\'errez, Hagenah,
      \cite{DiGuHa}) \label{theo;Gut}

        The existential theory with rational constraints of 
        a free group is decidable.
      \end{theo}

            This is why we wanted to define ``geometric'' elements in a local
            way. As we did, $\calL$ is a rational language. Indeed it is not
            hard to design a finite state automaton that recognize the
            presence of a sub-word of length at most $L$ that fails to define a
            $(\lambda,\mu)$-quasi-geodesic in $Cay \Gamma$, since there are
            only finitely many such possible words, and their list is
            computable with a solution to the word problem.

            From the Theorem above, and the equivalence of Proposition   
            \ref{prop;equiv}, we deduce:

            \begin{theo}
              The existential theory of a torsion free hyperbolic group (with
              parameters) is decidable.
            \end{theo}

         \paragraph{Looking further.}

            It appears, in fact, that such a result might not be sufficient to
            provide the main expected application: a significant
            simplification of Z.~Sela's algorithm on the Isomorphism Problem for
            torsion free hyperbolic groups. However, we have not demonstrated
            the full strength of the use of rational languages, but merely
            used the concept of geometric elements, and  $\calL_0$ among them.
            
            In this perspective, in \cite{DG} we use more refined languages, 
            in order to
            get more control (in particular inside 
            conjugacy classes of solutions). We  use the
            following statement, which is a direct consequence of the study above.

            \begin{prop} \label{prop;isom_hyp}

              Given $\G$ a torsion free hyperbolic group, generated by $S$ as
              above, there is an algorithm, explicit from the data of a presentation and the hyperbolicity constant, 
		that performs the following.

              Given a finite system $\mathscr{E}$ of equations in $\G$, and,
              for each  $\omega \in \Omega$, a rational language  
              $\mathcal{L}_{\omega}$  of geometric elements ({\it i.e. } 
              such that $\mathcal{L}_{\omega} \subset \calL$), the algorithm
              always terminates, and says ``yes'' if there exists a solution
              $s$
              of $\mathscr{E}$ with the property that, for all $\omega$,  
              \emph{any}
              pre-image of  $s(\omega)$ in $\calL$ is in
              $\mathcal{L}_{\omega}$, and only if  there exists a solution
              $s$ of $\mathscr{E}$ with the property that, for all $\omega$,  
              \emph{some}
              pre-image of  $s(\omega)$ in $\calL$ is in
              $\mathcal{L}_{\omega}$

            \end{prop}

		 The algorithm is explicit from the value of the hyperbolicity constant $\delta$ and the presentation of  $\G$. 
		Since there is an algorithm for computing $\delta$ from the presentation (see \cite{Papa}), we can actually deduce an algorithm that takes as input a presentation of a hyperbolic group, and a system of equations and inequations, and gives as output the answer to the satisfiability of the system.

            We will obtain such a statement in the generality of relatively
            hyperbolic groups, but for now, it is easy to see that a suitable 
            algorithm for Proposition \ref{prop;isom_hyp}
            is: find $\mathcal{E}'(c,rep)$ as above, and apply the 
            algorithm from \ref{theo;Gut} for the systems  $(\mathcal{E}'(c,rep), \mathscr{C}')$ where 
            $\mathscr{C}'= \{ \tilde{\omega} \in \mathcal{L}_{\omega}, \,
            \forall \omega \in \Omega \}$. 
		Answer yes if  one of these systems have a solution.

            We emphasize that 
            the algorithm can answer either yes or no when there is a solution 
            $s$ to $\mathscr{E}$ 
            with the property that, for all $\omega$,  
              \emph{some but not all}
              pre-image of  $s(\omega)$  in $\calL$ (not necessarily the canonical representative)  is in
              $\mathcal{L}_{\omega}$.
             Depending on the nature on the
            $\calL_{\omega}$, this uncertainty might never happen (this is the case for
            $\calL \setminus \calL_0$, since if a pre-image is inside, any
            pre-image is also inside), or happen only on finitely many
            solutions, or worse.

            It is also reasonable to ask whether such techniques can be used for
            the study of hyperbolic groups with torsion. Indeed, some specific
            cases can be handled, but the technical condition we need to
            require is rather unpleasant. The difficulty is the construction
            of canonical representatives. In what we develop, in the relative
            setting, the difficulty occurs precisely 
            in paragraph \ref{junctions}. The difficulty is already
            apparent for virtually free groups, although most of Rips and
            Sela's construction becomes trivial in this case: the associated 
            Bass-Serre tree almost provides canonical 
            representatives, but in fact
            only so-called canonical cylinders 
            (with ordered slice decomposition, see section 
            \ref{sec;cylsliced}), 
            from which it is sometimes
            impossible to extract a preferred path in the Cayley graph, in an
            equivariant way.

            Because it might shed light on the
            nature of the difficulty, let us try  to explain how far the
            construction can be generalized, but note that an argument of
            Delzant (ref  \cite[Rem. III.1]{Del}) shows that the 
            theorem of canonical
            representatives in its usual form cannot be extended to all
            hyperbolic groups.

            Let us qualify the action of $\G$ on $Cay \G$ as
            \emph{acylindrical}  if
            there is a number $\kappa$ such that the stabilizer
            of any finite subset of   $ Cay \G$ of diameter at least $\kappa$,
            is trivial.

            It is rather easy to adapt the construction of canonical
            representatives for hyperbolic groups with an acylindrical action
            of $Cay \G$, by considering, in paragraph \ref{junctions}, not
            slices, but finite unions of slices, so that stabilizers of
            these unions are trivial. Note also that Delzant's counterexample
            is virtually free, with a non-acylindrical action on its
            Bass-Serre tree.

            Although easy in the  hyperbolic case, this operation turns out to
            be rather delicate in the relative case.  The specific
            technical difficulties are due to so-called  angular slices, in the canonical cylinders. Nonetheless we claim that 
            it is still possible, but we chose 
            not to present the construction  here because, as far as we can
            see, it does not help to solve the case of torsion in generality,
            even in the hyperbolic case, 
            and is 
            perhaps of little interest, compared to its technicality.

   \section{Relative hyperbolicity, and paths in the coned-off graph}

      \subsection{Relative hyperbolicity, coned-off graph, angles and cones.}

        Let $\G$ be a finitely generated group, with a finite symmetric set of generators $S$. 
        Let $P_1, \dots, P_q$ be finitely generated subgroups of $\G$. We define, after Farb \cite{F}, 
        the \emph{coned-off} graph $\tilCaG$ of $\G$ with respect to the $P_i$ as follows. 
          Start with the Cayley graph $Cay \G$, 
          and  for each $k$, and each left coset $[gP_k]$ of $P_k$, add a vertex $v_{[gP_k]}$ together with 
          a family of edges $(v_{[gP_k]}, gh)$, for all $h\in P_k$. The graph obtained is denoted by  $\tilCaG$.

          Notice that the group $\G$ acts by left translations on the coned-off graph $\tilCaG$, 
          and for all $g$, the subgroup $gP_kg^{-1}$ fixes the vertex $v_{[gP_k]}$.  
          Whenever $P_k$ is infinite, the vertex  $v_{[P_k]}$ is of infinite
          valence (and also all of its translates), 
          and thus $\tilCaG$ is not locally finite.

          However, it may or may not be \emph{fine} (in Bowditch's terminology \cite{Brel}). 
          A graph is fine if any edge 
          is contained in only finitely many simple loops of given arbitrary length.

          We say that $\G$ is \emph{relatively hyperbolic} with respect to
          $P_1, \dots, P_q$, if $\tilCaG$ is hyperbolic and fine.

          This notion is indeed very robust, and rich. 
          We refer the reader to the now abundant literature on the subject for further discussion on properties, 
          naturality of the definition, and equivalent definitions 
          (see for instance, works of B.~Bowditch, B.~Farb, A.~Yaman,
          C.~Dru\c{t}u, M.~Sapir, D.~Osin, I. Bumagin, and myself). 

          In an arbitrary graph $G$ (but indeed we will be 
          interested in $\tilCaG$), 
          the \emph{angle} at a vertex $v$ between two
          edges $e_1=(v,v_1), e_2=(v,v_2) $
          adjacent at $v$ is the path distance (non-negative, possibly infinite) 
          in $G\setminus \{v\}$ between
          $v_1$ and $v_2$. It is denoted by $\Ang_v(e_1,e_2)$.

           Let now $\rho$ and $\theta$ be a positive number. 
           Given an oriented edge $e=(v,v')$ in $G$, the \emph{cone} centered
           at  $e$, 
           of radius $\rho$ and angle $\theta$, denoted by
           $\Cone_{\rho,\theta}(e)$ is 
           the set of vertices $w$ in $G$
          that  can be reached from $v'$ by a path  starting by  the edge $e$,
          of length at most $\rho$, and whose consecutive edges 
          make angles  at most $\theta$.

      \subsection{The group $\tilG$ labeling paths in $\tilCaG$.} \label{sec;tilG}

        In $\tilCaG$, there are two kinds of (oriented) edges. 
        Those that are in $Cay\G$ carry natural labels in $S$. 
        On the other hand,  those that are adjacent to a vertex  
        not in $Cay \G$ don't. 
        Nevertheless, if $e_1=(\gamma_1, v_{[gP_k]})$ and 
        $e_2 = (v_{[gP_k]},\gamma_2)$ are two adjacent such edges,  
        then we can associate to the ordered  pair $(e_1,e_2)$ a well 
        defined element of $P_k$, 
        namely $\gamma_1^{-1} \gamma_2 \in P_k$. We thus choose to label the
        (ordered) pairs 
        of consecutive edges containing a 
        same vertex $v_{[gP_k]}$ by elements of $P_k$.
        Clearly, this labeling is invariant by left translation.

        Let us now define the group $ \tilG = F_S * P_1 * P_2 * \dots * P_q$,
        where $F_S$ 
        is the free group over $S$. 
        For convenience of notations, we write $F_S = P_0$. 
        As $\tilG$ is a free product, every 
        $\tilde{\gamma} \in \tilG$  has a unique normal form 
        $ \tilde{\gamma} = x_1 \dots x_r$, where $\forall i, \exists j, x_i
        \in P_j \setminus\{1\}$, 
        and $x_i$ and $x_{i+1}$ 
          are in different factors, for all $i$. 

          Let us also define the
        \emph{long normal form} of   $\tilde{\gamma} \in \tilG$, to be the
        word obtained from its normal form after substituting to every $x_i
        \in P_0$, the normal form of $x_i \in P_0 = F_S$.

           There is a natural quotient $\tilG \tto \G$, 
          that is identity on $S$ and on each $P_i$. 
           We can interpret the elements of $\tilG$ as labeling paths
          in $\tilCaG$, via their long normal forms, according to 
          the labels defined above. Notice that paths of $\tilCaG$ starting at $1$, 
	and not using twice the same edge consecutively, are in one-to-one correspondence with long normal forms: 
	two different long normal forms (in the free product $\tilG$) will give two different paths and 
	any such path has a long normal form associated. 

          We thus denote by $\frap(\tilde{\gamma})$ the path in $\tilCaG$,
          starting at $1$, 
          and labeled by the long normal form of 
          $\tilde{\gamma}$.

          We define the \emph{sector} of angle
          $\theta$ of $P_k (<\tilG)$, denoted by $\Sec_k(\theta)$,
          to be the set of elements that label pairs of edges in
          $\tilCaG$, adjacent to $v_{[P_k]}$, with angle at most $\theta$ at $v_{[P_k]}$.

          It is easy to check that the fineness of $\tilCaG$  implies the
          finiteness of every cone in $\tilCaG$, and of every sector
          $Sec_k(\theta)$, for arbitrary $k$ and $\theta$.

          \subsection{Digression: speaking languages in groups.}

          Let us recall that a regular language on an alphabet is 
          one which is  the
          accepted language of a finite state automaton. 
          We refer to the first
          chapter of \cite{E} for a presentation of the theory.

          We recall definitions of a few class of languages in groups. 
           Let $G$ be a group, and $S$ be a finite symmetric 
           generating set, with
           an order. 

           A \emph{rational language} of $G$ is a subset of $G$ that 
           is the image
           of a regular language in the set of words  $S^*$ on 
           the alphabet $S$. 
           
           Given an element $g$ of $G$, the \emph{lexicographical normal form} 
           of $G$ is the word
           in $S^*$ that is the first in the lexicographical order, among the
           shortest words representing $g$. 
           A \emph{normalized regular language} of $G$ is a subset $A$ of $G$
           such that  the language  $\tilde{A}\subset S^*$, consisting of the
           lexicographical normal forms of the elements of $A$, is regular.

We mentioned normalized regular languages because in case $\tilG$ is a  free product of abelian groups, these languages are of rather easy use, and they allow a large number of possibilities, but the case we want to deal with is when $\tilG$ is more general.

Let us now explain another class of languages useful for general free products,
that was introduced by V.~Diekert and M.~Lohrey \cite{DiLo} for (even more general)
graph products of monoids. There are several simplifications in our
exposition, compared to \cite{DiLo}, due to the specific case we consider
({\it e.g. } for free products, the set they note $I$ is empty, their assumption
3.8 is always true for groups, their notations 
$U_\sigma$ and $V_\sigma$  are redundant with $M_\sigma$ for groups,
for a free product $IRR(R)$ is simply the set of normal forms, etc). Still we
introduce the next tools according to notations in \cite[Sec. 3.3]{DiLo}, to allow the comparison. 

Let $\Sigma$ be a finite set, and for all $\sigma \in \Sigma$, a group $M_\sigma$ is given. Let $\mathbb{F}(\Sigma)$ be the free product of the $M_\sigma$.

For each $\sigma$ let $\calC_\sigma$ be a family of  subsets of  $M_\sigma$, and  $\calC^\bullet_\sigma$ be the family of the same subsets where $1_\sigma$ has been removed in each of them. Let us write  $\calC^\bullet = \bigcup_\sigma  \calC^\bullet_\sigma$.

Let us define a $\calC^\bullet$-automaton to be a finite state automaton such
that every edge is labeled by an element of $\calC^\bullet$ ({\it i.e. } a
language in some $M_\sigma$). If $\mathbb{A}$ is such an automaton, the
language accepted by $\mathbb{A}$, noted $\calL(\mathbb{A})$  is defined as
usual, and is a subset of $(\bigcup M_\sigma^\bullet)^*$ where
$M_\sigma^\bullet = M_\sigma \setminus \{1_\sigma\}$: a word $a_1 \dots a_m \in
(\bigcup M_\sigma^\bullet)^*$ is accepted if there is a sequence of states
$q_0 \dots q_m$, with $q_0$ initial, and $q_m$ accepted, and transition edges
$e_i$ from $q_{i-1}$ to $q_{i}$, labeled by languages $L_i$ with $a_i\in L_i$.

The set of all accepted languages, for all possible  $\calC^\bullet$-automata is noted $\calL(\calC^\bullet)$, and is a subset of $\mathcal{P}( (\bigcup M_\sigma^\bullet)^*  )$. 

Note that, apart from languages accepted by $\calC^\bullet$-automata,
another remarkable element of  $\mathcal{P}( (\bigcup M_\sigma^\bullet)^*  )$
is the family  of  normal forms of elements of $\mathbb{F}(\Sigma)$, noted
$NF$, which is in natural bijection with $\mathbb{F}(\Sigma)$. Thus if $ L\in
\calL(\calC^\bullet)$,  the set $L \cap NF$ consists of all normal forms of
elements in $\mathbb{F}(\Sigma)$ that are accepted by the automaton.

Finally, we define $\mathcal{IL}(\calC) = \{ L \cap NF, \, L\in  \calL(\calC^\bullet) \}$. This is a family  of subsets of $(\bigcup M_\sigma^\bullet)^*$, and each such subset is included in $NF$, which allows to interpret it as a set of elements of  $\mathbb{F}(\Sigma)$.

Let us now consider our group $\tilG = P_0* \dots * P_q$, where $P_0$ is the 
free group
$F_S$. In the notation above, the
$M_\sigma$ are the factors $P_i$, $i\geq 1$, and the cyclic factors of $P_0$
generated by each element of $S$. 
 For each $i$, we choose
$\calC_i$ to be the family of finite or co-finite subsets of $P_i$. To justify this choice, we note that the decidability of the existential theory with parameters of $P_i$  is equivalent to the decidability of the existential theory of $P_i$ with constraints in $\calC_i$, which is the correct assumption for using the results in \cite{DiLo}.
 
 Let us denote by $\mathfrak{gnr}$ the union of the $\calC_i$.  The
 considerations above give rise to some specific class of languages, accepted
 by  $\mathfrak{gnr}^\bullet$-automata, denoted by 
$\mathcal{IL}(\mathfrak{gnr})$, and    let us then qualify these languages in $\tilG$ as
\emph{generalized normalized regular} (or $\mathfrak{gnr}$) languages. Let us emphasize that 
this class of languages is associated to a given free splitting of the group $\tilG$, and that
  a $\mathfrak{gnr}^\bullet$-automaton reads the long normal forms of elements of $\tilG$.

The aim of this tool is to define a class of constraints for which the
existential theory of $\tilG$ with constraints   would be decidable, provided that  the existential theory of each $P_i$, with parameters, is decidable. This will
follow from the main result of
 \cite{DiLo} (see \ref{theo;DL} below), and from the decidability of  the existential theory of $P_i$ with constraints in $\calC_i$ (which, as noted above, follows from our assumption, by choice of the $\calC_i$).

Let us try to clarify the difference between normalized regular
languages and  $\mathfrak{gnr}$ languages. The following lemma claims that any 
 $\mathfrak{gnr}$ language is normalized regular, but the converse is not
 true. In the case of free
 products of abelian groups, it can be interesting to use this latter
 class (together with  the main result of
 \cite{DiMu}, see \ref{theo;DM} below) , since it gives more latitude in the
 choice of constraints. On the other hand, when the free product has non-abelian factors, 
 normalized regular languages may not be  useful since the main theorem of
 \cite{DiLo} may not apply.

 \begin{lemma} With the above assumptions, given arbitrary orders in
   generating sets  of each $P_i$,  any $\mathfrak{gnr}$ language is normalized regular. 
 \end{lemma} 

 Consider a language $L$ in  $\mathcal{IL}(\mathfrak{gnr})$, and an associated
 $\mathfrak{gnr}^\bullet$-automaton $\mathbb{A}$. Given an (oriented) edge $e$
 of the automaton, it is associated to a finite or co-finite subset $X$
 of some $P_i$. 
 In both cases, there is an automaton  $A_e$ whose alphabet is the generating set of $P_i$, and whose accepted language is that of lexicographical normal forms of elements of $X$.

Replace then the edge $e$ by an automaton $A_e$ in the following way:   the initial state of $A_e$ is identified to  the origin state of $e$,  and the end state of $e$ is duplicated and identified to each accepted state of $A_e$, all other state of $A_e$ being non-accepted. Thus, it requires a lexicographical normal form of an element of $X$ to go through this automaton. 

Perform the construction for each edge of the original
$\mathfrak{gnr}^\bullet$-automaton $\mathbb{A}$, the result is a finite state
automaton whose  alphabet is the union of generating set of the $P_i$.

  It is then  a direct consequence of the construction that  a lexicographical normal form is accepted if and only if  it defines an element in $L$. Thus, $L$ is normalized regular. $\square$

          \subsection{Choosing a language without detour.}

            The remaining of this section is devoted to the choice of a good
             class of ``geometric'' elements in $\tilG$ (similarly to the
             hyperbolic case in the previous section). We look for a class of
             elements that label geometrically relevant paths in $\tilCaG$ and
             that is recognizable by a finite state automaton, in the sense of
             generalized normalized regular languages.

            \subsubsection{Detours and small detours.} 

            Let now 
            $\tilde{\gamma}\in \tilG$, with normal form $\tilde{\gamma}
            = x_1 x_2 \dots x_r$. 
            Let $i_0\geq 2$ and $s \leq r-i_0 -1$, and let
            $\tilde{d} = x_{i_0} x_{i_0 +1} \dots x_{i_0 +s}$. 

            We say that
            $\tilde{d}$ is a $\theta$-\emph{detour} of 
            $\tilde{\gamma}$ if there is
            $k$ such that 
            both $x_{i_0-1}$ and $x_{i_0+s+1}$ are in $ P_k \setminus
            \Sec_k(\theta)$ and if  $\frap(\tilde{d})$ starts and ends at
            distance $1$ from $v_{[P_k]}$ in $\tilCaG$. We call length of 
            a detour, the length of the path it labels in $\tilCaG$.
            Intuitively, if  $\tilde{\gamma}$  has a detour, the path $\frap(\tilde{\gamma})$ contains a severe
            backtracking on some vertex translated of $v_{[P_k]}$.
            We also say that a path 
            has a $\theta$-detour if it
            is labeled by an element of $\tilG$ with a  $\theta$-detour.

            Let $\mu \geq 1$. Let us consider  a $\theta$-detour $\tilde{d}$ of
            $\tilde{\gamma}$, with the notations as above. We say that it is
            $\mu$-\emph{small} if all the $x_{i_0+j}$, $j=0, \dots, s$ are in
            $F_S \cup (\displaystyle \bigcup_{k=1}^q
            \Sec_k(5\mu(\mu+\theta)))$.  Intuitively, a small detour is a detour
            without large angle inside.

            The choice of constants  will be justified
            in the following.

            \subsubsection{The language $\calL$.}

            We want to define a class of elements of $\tilG$, that we wish to
            call geometric, in the sense that they are relevant to the
            geometry of $\tilCaG$.

            Let  $M$ be a bound on the cardinality of 
            cones of radius and angle $50\delta$. 
            Let $L_1 = 10^4 \delta M$, $L_2= 10^6 \delta^2 M$, and let $\theta$
            be the computable constant that we denoted by $\Theta$ in \cite{Dis} (for instance
            $10^4(D+60\delta)$ is suitable, where $D$ is a fellow traveling constant for
            $1000\delta$-quasi-geodesics, and greater than any angles at finite
            valence vertices).  
           
            As already mentioned, it is a classical property of  
            hyperbolic geometry that there
            exists\footnote{In fact there are explicit 
              values in \cite[p.30-31]{CDP}} $L$, $L_1'>0$ and $L_2'$ 
            such that any $L$-local $(L_1,L_2)$-quasi-geodesic is
            a (global)  $(L_1',L_2')$-quasi-geodesic. One can
            choose $L>L_2'$.

            We then qualify as \emph{geometric}, the elements $\tilde{\gamma}
            \in \tilG$ that do not have any  $\theta$-detour   
            such that $\frap( \tilde{\gamma})$ is a $L$-local
            $(L_1,L_2)$-quasi-geodesic. We denote by $\mathcal{L} \subset
            \tilG$ the set of geometric elements of $\tilG$.

            The aim of the next paragraphs is to prove that 
            $\mathcal{L}$ is generalized normalized regular, 
            and that so
            is the subset of $\mathcal{L}$ representing the 
            trivial element of $\Gamma$. 

            It will turn out that there are good canonical representatives in
            $\calL$ for turning equations in $\G$ 
            into equations in $\tilG$, hence the interest of this language.

            \subsubsection{A local characterization of $\calL$.}

            The aim of
            this part is to prove the following characterization of $\calL$, that
            weakens the absence of detours to the absence of small detours.

            \begin{prop}\label{prop;caracL}
              An element $\tilde{\gamma}$ is in $\calL$ if and only if it has
              no $L_2'$-small $\theta$-detour of length at most $L_2'$, 
              and $\frap( \tilde{\gamma})$ is a $L$-local
              $(L_1,L_2)$-quasi-geodesic.
            \end{prop}
          
            This characterization will be used in order to prove the regularity
            of $\calL$, whereas we will need the property stated as definition in
            order to show later the fellow traveling property of 
            elements of $\calL$.  

            Before proving the characterization itself, 
            we establish a couple of lemmata.

        \begin{lemma}(Loops in absence of detours)\label{lem;loop}

          Let $\theta>0$, and consider the family of 
          oriented edges of 
          an oriented path of length $T$, without $\theta$-detour.

          Let $v$ be a non-extremal vertex on the path, 
          and let us denote by $e_1, \dots,
          e_k$ the edges of our family that arrive at $v$, and $f_1,  \dots
          f_k$ be the edges of our family that start at $v$.
  
          Let $a >0$. If $\Ang_v(e_1,f_k) \leq a$, then  
          $\Ang_v(e_1,f_1) \leq a+T(T+\theta)$. 
          Moreover for all $i\leq k-1$,  $\Ang_v( f_i, e_{i+1}) \leq T$ for
        all $i\leq k-1$.

        \end{lemma}

        {\it Proof. } 
        First let us notice that $f_i$ and $e_{i+1}$ are the first and last
        edges of a loop that passes only once at $v$ and that has length at
        most $T$. This provides immediately $\Ang_v( f_i, e_{i+1}) \leq T$ for
        all $i\leq k-1$.

          By triangular inequality for angles, 
        $ \Ang_v(e_1,f_1) \leq  \Ang_v(e_1,f_k)  + \displaystyle
        \sum_{i=1}^{k-1}    \Ang_v( f_i, f_{i+1})      \leq   \Ang_v(e_1,f_k)  + \displaystyle
        \sum_{i=1}^{k-1} ( \Ang_v( f_i, e_{i+1}) + \Ang_v( e_{i+1}, f_{i+1}))
        $, and therefore  $ \Ang_v(e_1,f_1) \leq  a +(k-1)T  + \displaystyle
        \sum_{i=1}^{k-1} \Ang_v( e_{i+1}, f_{i+1}))$.

        We can assume that $\Ang_v(e_1,f_1)>\theta$, otherwise there is
        nothing to prove. Then if for some $i\leq k-1$, $\Ang_v( e_{i+1},
        f_{i+1})) >\theta$, this produces a $\theta$-detour for our path,
        between the edges $f_1$ and $e_{i+1}$, which is opposed to 
        our assumption. Thus, for all $i$, $\Ang_v( e_{i+1},
        f_{i+1})) \leq \theta$, and therefore, $\Ang_v(e_1,f_1) \leq  a
        +(k-1)(\theta +T)$. Since $k\leq T$, this proves the lemma. $\square$

    \begin{lemma} \label{lem;findsmall}
            Let $\tilde{\gamma} \in \tilG$,  $\theta>0$, and $\mu>0$. 

            Any
            $\theta$-detour $\tilde{d}$ of $\tilde{\gamma}$, such that the
            length of $\frap(\tilde{d})$ is at most $\mu$,
            contains, or is, a $\mu$-small $\theta$-detour.

          \end{lemma}

          {\it Proof. } Let us argue by contradiction, and 
          consider $\tilde{d}= x_{i_0} x_{i_0 +1} \dots x_{i_0 +s}$ 
          a $\theta$-detour of length at most $\mu$, without $\mu$-small $\theta$-detour. Let
          us assume that $\tilde{d}$ has minimal length $s$ among
          $\theta$-detours of $\tilde{\gamma}$. Since it is assumed not to be small, 
          one letter $x_{i_0+j}$ is in some $P_{k'} \setminus \Sec_{k'}(5\mu(\mu+\theta))$.
 
          By definition of detour, the path $\frap(\tilde{d})$, together with
          two edges from its end
          points to $v_{[P_k]}$,  defines a loop, which we denote by
          $\frap(\tilde{d})^+$,  
          of length at most $\mu+2$. 
          Because of $x_{i_0+j}$, this loop  contains a vertex
          $v$ at which its two consecutive edges make an angle of at least
          $5\mu(\mu+\theta)$.

          Since this angle is greater than $\mu$,  the end points of the edges
          cannot be joined by a loop of length $\mu$ avoiding $v$, hence, 
          $\frap(\tilde{d})^+$   has  to pass at least twice on $v$.
 
          One can then apply Lemma \ref{lem;loop} to $\tilde{d}$, for which
          $T\leq 2\mu$, at the vertex $v$. One can use $a= \mu$, since the
          first and last edges  of $\frap(\tilde{d})$ containing $v$ are in a
          sub-loop of  $\frap(\tilde{d})^+$ of length at most $\mu$ that does
          not pass twice by $v$, and make therefore an angle at most $\mu$. 
          One thus gets
          $ 5\mu(\mu+\theta) < \mu + 2\mu(2\mu+\theta)$ which is a
          contradiction. $\square$

          Let us now prove Proposition \ref{prop;caracL}. By definition an
          element $\tilde{\gamma}$ is in $\calL$ if and only if it has no
          $\theta$-detour, and if $\frap(\tilde{\gamma})$ is a $L$-local 
          $(L_1,L_2)$-quasi-geodesic. Any such element satisfies the
          condition of  Proposition \ref{prop;caracL}, 
          and for the converse, given
          an element such that  $\frap(\tilde{\gamma})$ is a $L$-local 
          $(L_1,L_2)$-quasi-geodesic, since $\frap(\tilde{\gamma})$ is a
          global $(L_1',L_2')$-quasi-geodesic, any detour is of length at
          most $L_2'$, and therefore, by Lemma \ref{lem;findsmall} 
          (with $\mu=L_2'$), any
          $\theta$-detour contains a $L_2'$-small $\theta$-detour, of length
          at most $L_2'$. Hence the equivalence stated. $\square$

          \subsubsection{Normalized regularity of $\calL$.}

        \begin{prop} \label{prop;reg}
          The language $\mathcal{L} \subset \tilG$ is generalized normalized regular 
          (hence in particular normalized regular).
        \end{prop}
        
        {\it Proof.}
        We need an automaton that reads the long normal forms of elements of
        $\tilG$ and accepts exactly  those in $\calL$.

        Let $\mathscr{D}$ the list of all long normal forms $w$ of 
        elements of $\tilG$ 
        that are $L_2'$-small $\theta$-detours 
        with $\length(\frap(w))\leq  L_2'$ 
        (as in Proposition \ref{prop;caracL}). 
         The long normal form of a small detour
         consists of elements of $S \cup (\displaystyle \bigcup_{k=1}^q
          \Sec_k(5L_2'(L_2'+\theta))$. 
          Therefore $\mathscr{D}$ is finite and computable (indeed the sets  $\Sec_k(\kappa)$ are computable if one knows the list of simple loops of $ \widehat{Cay}\G$ of length $2\kappa$, which we know, or are given, by assumption).

        It is then easy to design a $\mathfrak{gnr}^\bullet$-automaton  
        that checks
        whether a long  
        normal form of an element contains a
        $L_2'$-small detour of length at most $L_2'$. Indeed it suffices to
        check whether the given word contains a sub-word in
        $\mathscr{D}$, preceded and followed by two elements in some
        same $P_k\setminus \Sec_k(\theta)$ (which is co-finite in $P_k$).

        It remains to design a $\mathfrak{gnr}^\bullet$-automaton that, given  a
        long normal form without $L_2'$-small $\theta$-detour of
        length less than $L_2'$, checks whether it labels a $L$-local
        $(L_1,L_2)$-quasi-geodesic in $\tilCaG$.

        A  path 
        fails to be a  $L$-local $(L_1,L_2)$-quasi-geodesic 
        if and only if one of its sub-paths of length  $\ell \leq L$ is labeled by an element
        $g_L \in \tilG$ 
        satisfying in $\tilCaG$ the inequality  
        $d(1,\overline{g_L}) < (\ell -L_2)/L_1$. 
        We say that the long normal form of such a $g_L$, 
        is \emph{forbidden}.

        Moreover, since $L>L'_2$, a path of length at 
        most $L$ without $\theta$-detour have no $L_2'$-small 
        $\theta$-detour of
        length less than $L_2'$. 

        These observations reduce the problem to finding an automaton 
        that checks whether  a
        long normal form $w$ with  $\length(\frap(w))   \leq L$, 
        and without $\theta$-detour, is forbidden.

        Let $A= 2(L+\theta)^2  $. The  long normal form $w$ of an 
        element of $\tilG$ 
        can be written in a unique way
        $w= w_1a_1w_2 \dots a_{k-1} w_k a_k w_{k+1}$, for some $k$,  
        with, for all $i$,  
        $a_i$ an element in $P_k\setminus \Sec_k(A)$ for 
        some $k\leq q$ ({\it i.e. } 
        defining in $\widehat{Cay}\Gamma$ an 
        angle greater than    
        $A$)      
        maximal for this property, and  with $w_i$ the 
        long normal form of an element defining a path 
        in $\widehat{Cay}\Gamma$ with no   angle greater than    
        $A $.  By maximality of the $a_i$, one has 
        $\length(\frap(w))= 2k+\sum length(\frap(w_i))$.

        Let us now prove the following lemma, before finishing the proof of
        Proposition \ref{prop;reg}

        \begin{lemma} Assume that the element defined by the long  normal
        form $w$ has no $\theta$-detour, and that $\length(\frap(w))\leq L$. 
        Then, with the notations as above, 
        it is forbidden if and only if 
       $(2k+\sum length(\frap(w_i)) -L_2)/L_1> 2k+\sum d(1,\bar{w_i})$,  
       where $\bar{w_i}$ is the element of $\Gamma$ defined by $w_i$.  
       \end{lemma}

       The word $w$ labels a path $\frap(w)$ in $\tilCaG$ of 
       length $\length(\frap(w)) = 2k+\sum
       \length(\frap(w_i))  $
       from $1$ to $\bar{w}\in \G$. 
       It is forbidden if and only
       if $\length(\frap(w)) > L_1 d(1,\bar{w}) + L_2 $. 
       It remains to show that $d(1,\bar{w}) \geq 2k+\sum d(1,\bar{w_i})$.

       Let us denote by $v_i$, $i=1,\dots, k-1$,  the
       consecutive vertices of $\frap(w)$ where the letters $a_i$ label an
       angle. In $\tilCaG$, for $i\in [1,  k-1]$,  the word $w_i$ labels a path between finite valence vertices 
       adjacent to $v_i$ and $v_{i+1}$, thus  $d(v_i,v_{i+1}) = d(1,\bar{w_{i+1}})+2$ for all $i\in [1,  k-1]$. 
       Similarly, 
       with the convention $v_0=1$ and $v_{k+1} =
       \bar{w}$, one has $d(v_0,v_{1}) = d(1,\bar{w_1})+1$, and $d(v_k,v_{k+1}) = d(1,\bar{w_{k+1}})+1$.

         Thus, $\sum_{i=0}^k  d(v_i,v_{i+1})  = 2k+\sum
       d(1,\bar{w_i})$. Therefore it suffices to check that $ \sum_{i=0}^k  d(v_i,v_{i+1})  = d(1,\bar{w})$, 
       which is true if 
       all the $v_i$ are on
       some geodesic $[1,\bar{w}]$. 

       Consider a geodesic $[1,\bar{w}]$, and
       assume that  $v_i \in [1,\bar{w}]$ for some $i$.
       The path defined by the union of the sub-path $\frap(w)$ from
       $1$ to $\bar{w}$, and of 
       a geodesic $[1,\bar{w}]$, is
       a loop of length at most $2L$.  
       Since the angle of  $\frap(w)$ at
       $v_{i}$ is greater than $A>2L$, 
       the loop has to pass at least twice
       on this point.
       Since
       $v_{i} \notin [1, \bar{w}]$, the path $\frap(w)$ passes twice at
       $v_{i}$. One can then apply Lemma \ref{lem;loop} for $a=2L$ to obtain
       that the angle at $a_{i}$ is at most $2L+L(L+\theta)$. 
       Since it is also
       at least $A=2(L+\theta)^2 $,  this 
       contradicts  the choice of $A$. $\square$

       Let us then finish the proof of Proposition \ref{prop;reg}. Recall that
       it is enough to find an automaton that,  given a  long normal
       form  $w$ with $\length(\frap(w))$ at most $L$, 
        without $\theta$-detour, checks whether it is forbidden.

       Consider the set $\mathcal{S}$ of all the words in long normal form, that define in 
       $\widehat{Cay}\Gamma$ a path of length at most $L$ and of maximal angle at most $A$. 
       This set is finite, 
       and for all $k \leq  L/2 $, the set $\mathcal{S}_k$ of $(k+1)$-tuple 
       of elements of $\mathcal{S}$ that 
       satisfy the inequalities  $(2k+\sum length(\frap(w_i))-L_2)/L_1> 2k+\sum d(1, \gamma_{w_i})$ 
       and $2k+\sum length(\frap(w_i))\leq L$  is finite and computable. The lemma asserts that      
       the     forbidden words without $\theta$-detour  
       are exactly  the concatenations the elements $w_i$ of a certain set 
       $\mathcal{S}_k, k\leq L/2$ with  
       elements of the factor groups  in one of the $P_k\setminus \Sec_k(A)$ (which are co-finite in the $P_k$). 
       This can be easily recognized
       by a $\mathfrak{gnr}^\bullet$-automaton. $\square$

          \subsection{Conical fellow-traveling, 
            and finiteness of $\calL_0$.}

          We now prove  that
          quasi-geodesics without small detours fellow travel in a strong
          sense.

          \begin{prop} \label{prop;stab_con}

            Let $\lambda, \mu, \chi$ be positive constants, and let
            $\epsilon$ be a constant of fellow traveling for
            $(\lambda,\mu)$-quasi-geodesics in $\widehat{Cay}\Gamma$, which is
            hyperbolic. 
            Set $\epsilon'=   4(\epsilon +   \lambda(3\epsilon + \mu)+   
            ( \lambda(3\epsilon +\mu) )  (\lambda(3\epsilon +\mu )+  \chi))
            +50\delta $.

            Let $\tilde{\gamma}$ be an element of $\tilG$, without 
            $\chi$-detour.
            Assume moreover that $\frap(\tilde{\gamma})$ is a
            $(\lambda,\mu)$-quasi-geodesic of $\tilCaG$. 
            Denote by $1$ and $v$ its starting
            and ending point, and assume $v\neq 1$.

            Then, if $[1,v]$ is a geodesic segment with set of edges 
            $\mathcal{E}$, the path
            $\frap(\tilde{\gamma})$ is contained in the union 
            $\bigcup_{e\in \mathcal{E}} \Cone_{\epsilon' ,\epsilon'} (e)$.

          \end{prop}

        Let $\rho:[0,T] \to \tilCaG$ be an arc-length parameterization of  $\frap(\tilde{\gamma})$.  
        Let $v= \rho(t)$, and assume it is not in  $[x,y]$.  
        Let us define $t_- = \max\{0, t - \lambda(3\epsilon +\mu )\}$, and $t_+ = \min\{ T, t+\lambda (3\epsilon + \mu ) \}$.
        Let us also choose $\tau_-$ a geodesic from $\rho(t_-)$ to a closest point $x'$ of $[x,y]$, 
        and $\tau_+$ in the same way. The choice of $t_-$ and $t_+$ guarantees that $\tau_-$ and $\tau_+$ 
        are disjoint, and  that $[x',y']$ is not reduced to a point. 
        We now choose $t'_-$ such that $\rho(t'_-)$ is the last common point of $\rho$ and $\tau_-$, and $\tau'_-$ is the subsegment of $\tau_-$ joining $\rho(t'_-)$ to $x'$ (and similarly for $t_+$). Since the length of $\tau_-$ and $\tau_+$ is at most $\epsilon$, one still has $t \in [t'_-,t'_+]$.  
        We now consider the loop  equal to the concatenation $\tau'_-, \rho|_{[t'_-,t'_+]}, \tau'_+, [y',x']$. It contains an edge $e$ of $[x,y]$ and the vertex $v$, and it is of length at most $2\epsilon + 2  \lambda(3\epsilon +\mu ) + d(x',y')$, that can be bounded by  $4\epsilon + 4  \lambda (3\epsilon + \mu )$. 
        
        We now need to bound the angles in this (not necessarily simple) loop. 
        The multiple points of the loop are all on $ \rho|_{[t'_-,t'_+]}$. 
        This path has no $\chi$-detour, and is of length at most   $2   \lambda(3\epsilon +\mu)$,   therefore, 
                        by Lemma \ref{lem;loop}, its angles on multiple points
        are at most $   4\epsilon + 4  \lambda(3\epsilon +\mu)+   (2
        \lambda(3\epsilon +\mu) )  (2   \lambda(3\epsilon +\mu)+     \chi)
        $. On the other hand the  angle on every simple point of the loop  is
        at most   $   4\epsilon + 4  \lambda (3\epsilon +\mu )$. This provides
        a path from the edge $e\in [x,y]$ to $v$ of length at most
        $4\epsilon + 4  \lambda(3\epsilon +\mu)$ and maximal angle at most 
        $4\epsilon + 4 \lambda(3\epsilon +\mu)+   (2  \lambda(3\epsilon +\mu) ) 
        (2 \lambda (3\epsilon +\mu )+ \chi) \leq\epsilon'$. 
        Therefore, $v \in Cone_{\epsilon',\epsilon'}(e)$.   $\square$

        \begin{coro}  \label{coro;L_0}
          The language $\mathcal{L}_0$ of the elements of $\mathcal{L}$
          that represent the identity of $\Gamma$ is finite and computable, and $\calL
          \setminus \calL_0$ is generalized normalized regular.
        \end{coro}

       {\it Proof. } For the choice of
       constants $\lambda = L_1, \mu= L_2, \chi = \theta$, 
       Proposition \ref{prop;stab_con} indicates
       that for all $\tilde{\gamma} \in \calL_0$, the path $\frap(\tilde{\gamma})$ 
       remains in  $\Cone_{\epsilon',\epsilon'}(e)$, where $e$ is
       the first edge  of $\frap(\tilde{\gamma})$, being in particular one of the finitely
       many edges adjacent to the base vertex $1$ of $\tilCaG$.   

       Since it is a $(L'_1,L'_2)$-quasi-geodesic, and a closed loop, 
       it has length at most $L_2'$. This leaves only finitely
       many possibilities for $\frap(\tilde{\gamma})$, therefore for
       $\tilde{\gamma}$. 
       The rest of the corollary 
       follows immediately from Proposition \ref{prop;reg}. $\square$

   \section{Canonical Representatives in $\tilG$}
   
      We are now interested in providing canonical representatives, in the way
      of Rips and Sela \cite{RS}, 
      for elements in $\G$. We want these representatives to be ``geometric'',
      that is to belong to $\calL$, in order to have control on them.
      
      Fortunately, the most technical part was already done in our previous
      work \cite{Dis}. We recall the main result we need.

      \subsection{Sliced cylinders.} \label{sec;cylsliced}

   Let us briefly recall  the results of our constructions of \cite{Dis} 
      (where details can be found).
      In the following, for $a$ and $b$ in 
      $\widehat{Cay} \Gamma$, and $l \in \mathbb{N}$, the set  $Cyl_l(a,b)$ 
      (read ``cylinder of parameter $l$ for the pair $(a,b)$'')  
      is a finite, computable subset of 
      $\widehat{Cay} \Gamma$, lying in a union of cones of known 
      radius and angle, and 
      centered on the edges of a geodesic path from $a$ to $b$.
      
      Any cylinder has an ordered decomposition in slices, that we write
      $Cyl_l(a,b)= \bigsqcup_{i=1}^s S_l(a,b)_i$   \cite[Def. 2.16]{Dis}.  
  
      These slices are  equivalence classes for a certain 
      relation (defined by a cocycle  \cite[Def. 2.15]{Dis}). 
      There exists a universal constant  such that any two
      consecutive slices of a cylinder are in a cone of that angle and radius  \cite[Corollary 2.20]{Dis}. 

      In particular, in $\tilCaG$,  there are only finitely many orbits of subsets that are
      slices of some cylinder, and also of subsets that are pairs of
      consecutive slices in some cylinder.

      The definitions of these slices and cylinders are compatible with
      isometries, so that one has the important property that $Cyl_l(a,b) =
      \gamma^{-1} Cyl_l(\gamma a, \gamma b)$ for all $\gamma\in \Gamma$ (see
      \cite[Lemma 2.7]{Dis}).

      The main property of cylinders and slices 
      is a stability property, which we rephrase now,  
      that states that they almost coincide in selected triangles.

      \begin{theo} (\cite[Theorems 2.9 and 2.22]{Dis}) \label{theo;cancyl}

        Let $\Gamma$ be a  relatively hyperbolic group, with a finite generating set. 
	Let $\delta$ be a
      hyperbolicity constant of a coned-graph. 
  	There exists computable constants $\kappa_1, \kappa_2$, depending only on
      the generating set such that the following holds.

        Let $\mathcal{A}$ be an arbitrary finite subset of $\G$, $n$ be an integer, and 
        consider $n$ words of three letters in  
        $\mathcal{A}$, $\tte_i= a(i)_1\, a(i)_2\, a(i)_3$,  with $i\leq n$,
        such that, for all $i$, 
        $\tte_i  \underset{\G}{=}   1$ in   $\Gamma$.

        Then, there exists an integer $l \in [1,
        10^{30}(\delta+1)^{10}]$, such that, 
        for all $i\leq n$,  the unique ordered slice decomposition of
        $Cyl_l(1,a(i)_j)$ can be written  

        $$ Cyl_l(1,a(i)_j)=
        (\bigsqcup_{k=1}^{l(j)} L(j)_k) \sqcup ( \bigsqcup_{k=1}^{c(j)}
        C(j)_k  )\sqcup  ( \bigsqcup_{k=1}^{r(j)} R(j)_k  ) $$
   
        where the following equalities hold: $r(j)=l(j+1 \, [3])$ for
        $j=1,2,3$, and
        $(a(i)_j)^{-1} R(j)_k = L(j+1\, [3])_{r(j)-k+1}$ for all $j\leq 3$, and all $k$,
        and where moreover, $c(j)\leq \kappa_1 n + \kappa _2$ and no 
        slice $C(j)_k$ is an
        angular slice of angle greater than $\kappa_2$. 
        
      \end{theo}

       In other words, the theorem states that for the considered relations
       $\tte_i= a(i)_1 \, a(i)_2\, a(i)_3$, there is $l$ such that the cylinders
       are of the form:
 
       $$
       \begin{array}{lllllllllll}
         Cyl(1,a(i)_1) & = &
         (S_1, & S_2,  & \dots, & S_k, &\, \mathcal{H}_z,\, & 
         T_m,  & T_{m-1}, & \dots, & T_1)   \\
         Cyl(a(i)_1,a(i)_1a(i)_2) & = &
         (S_1, & S_2,  & \dots, & S_k, &\, \mathcal{H}_y,\, & 
         V_p, & V_{p-1}, & \dots, & V_1)   \\
         Cyl(a(i)_1a(i)_2,1) & = &
         (T_1, & T_2,  & \dots, & T_m, &\, \mathcal{H}_x,\, & 
         V_p, & V_{p-1}, & \dots, & V_1),
       \end{array} 
       $$

        where $S_1, \dots, S_k, T_1,\dots, T_m$ and  
        $V_1,\dots, V_p$ are slices and where each $\mathcal{H}_v$,
        $(v=x,y,z)$ is a set of at most $\kappa_1 n + \kappa _2$ 
        consecutive slices,
        without parabolic slice of angle more than $\kappa_2$. 

        The constants
        $\kappa_1$ and $\kappa_2$ can easily be made explicit from \cite[Theorems 2.9, 2.22]{Dis}, 
	to which we refer the reader (note that they do not depend on the parameter noted $l$ in \cite[Theorem 2.9]{Dis}, 
	which may not be computable).

      \subsection{Centers and junctions between slices.}\label{junctions}

      Let $\mathscr{F}_1$ be a finite family of  
      representative of $\Gamma$-orbits of subsets of 
      $\widehat{Cay} \Gamma$ that are slices of some cylinder.

      Let us make the choice, for each $\sigma \in \mathscr{F}_1$, 
      of a point 
      $x_\sigma  \in \sigma$, which we call the {\it center} of $\sigma$. We
      can choose it to be a point of $\G$ as soon as $\sigma$ is not reduced
      to a single parabolic point, since in this case, a slice always contains
      a non-parabolic point.

      Similarly let $\mathscr{F}_2$ be a finite 
      set of representatives of equivalence classes of $\Gamma$-orbits of
      pairs of points belonging to two consecutive slices in some
      cylinder. The existence of such a finite set is ensured by the fact that
      any slice is finite, and that there are finitely many orbits of pairs of
      consecutive slices (as already mentioned).

      For all pair of points $x_1, x_2$ in  $\mathscr{F}_2$,  let us 
      choose a geodesic path $\mathfrak{j}_{(x_1,x_2)}$ in $\widehat{Cay} \Gamma$ 
      joining them,  
      and let us call this path a {\it junction}.  

      Note that by 
      \cite[Corollary 2.20]{Dis}, every junction is shorter than $1000\delta$, and has maximal
      angle at most $\theta$.

      \begin{lemma} \label{lem;junctions}
        Let $\G$ be a torsion free relatively hyperbolic group. 
        Let $Cyl_l(x,y)$ be a cylinder in $\tilCaG$ with ordered 
        slice decomposition $Cyl_l(x,y) = \sqcup_{j=1}^s S_j$.

        For all $j\leq s-1$, there is a unique pair of points $(x_j, x_{j+1}) \in S_j \times
        S_{j+1}$, such that there is $\gamma_j, \gamma_{j+1}$ such that $\gamma_j S_j \in \mathscr{F}_1$ has center $\gamma_j x_j$, 
        and similarly for $j+1$.
     
        There is also  a unique segment $\mathfrak{j}_{(x_j,x_{j+1})}$,
        depending only on $S_j$ and $S_{j+1}$, 
        such that there is $\gamma \in \G$ with $\gamma (x_j, x_{j+1}) \in \mathscr{F}_2$
        and  
        $\gamma \mathfrak{j}_{(x_j,x_{j+1})} = \mathfrak{j}_{(\gamma x_j, \gamma
        x_{j+1})}$.

      \end{lemma}

      {\it Proof. } 
      The existence of such objects is clear since  $\mathscr{F}_1$ and
      $\mathscr{F}_2$ are orbit representatives. 
      The interesting issue is uniqueness. Assume that  $x_j', x_{j+1}'$ and $
      \mathfrak{j}'_{(x_j,x_{j+1})} $ also satisfy the conclusion. 

      Let us first prove that $x_j = x'_j$ (and the same will hold for
      $x'_{j+1}, x_{j+1}$). By assumption there are $\gamma$ and
      $\gamma'$ such that $\gamma S_j=\sigma \in \mathscr{F}_1$, 
      of center $\gamma x_j$, 
      and $\gamma' S_j =\sigma $, of center $\gamma'
      x'_j$. Therefore, $\gamma'^{-1} \gamma \sigma =\sigma$, and
      $<\gamma'^{-1} \gamma>$ has a finite orbit in $\tilCaG$. This implies
      that either $ \gamma'^{-1} \gamma$ has finite order or is parabolic. 
      If it has finite order, by absence of torsion, $ \gamma' = \gamma$, and  $x_j = x'_j$. 
      If it is parabolic (of infinite order), it fixes a single point in $\tilCaG$, and it is its only finite orbit. 
      Thus $\sigma$ is reduced to this single point, and clearly  $x_j = x'_j$.

      Now, there is a unique pair in  $\mathscr{F}_2$
      in the orbit of $(x_j,x'_j)$, by definition of $\mathscr{F}_2$. 

      Therefore, there are
      $\gamma$ and $\gamma'$ in $\G$ such that $\gamma
      \mathfrak{j}_{(x_j,x_{j+1})} =  \gamma' \mathfrak{j}'_{(x_j,x_{j+1})}$.
      This means that $\gamma'^{-1} \gamma $ fixes two different points in
      $\tilCaG$, and therefore that it is in two distinct parabolic
      subgroups. This easily implies that it has finite order, and 
      since $\G$ is torsion free, $\gamma = \gamma'$ hence  
      $\mathfrak{j}_{(x_j,x_{j+1})}= \mathfrak{j}'_{(x_j,x_{j+1})}$.
        $\square$

      With the notations of the lemma above, we call {\it junction} 
      the geodesic      path  $\mathfrak{j}_{(x_j,x_{j+1})}$.  

   \subsection{Canonical representatives.}

      There are two possibilities for a junction: either it is between two
      points of $\G$, and it is labeled by a well defined element of $\tilG$,
      or one of its ends is a parabolic point. Note that it is not possible
      that both its ends are parabolic points since in the consecutive slices of a cylinder there are no
      consecutive parabolic slices (see \cite[Def. 2.6]{Dis}).
      In this case, the two consecutive junctions at the parabolic point
      provide a path from two points of  $\G$, and it corresponds to a 
      well defined element of $\tilG$.
      
      This allows to define 
      $\mathfrak{rep}_l(\gamma) \in \tilG $
      the canonical representative for $\gamma \in \G$ (for the
      parameter $l$) to be the product in $\tilG$ of all the elements
      associated to junctions or pairs of consecutive junctions between the
      ordered slices of $Cyl_l(1,\gamma)$.
      Theorem \ref{theo;cancyl} therefore gives  the following.

      \begin{theo}  \label{theo;repcan}

        Let $\Gamma$ be a torsion free  relatively hyperbolic group, with notations as above. There
        exists an explicit constant $\kappa$ such that the following holds.

        Let $\mathcal{A}$ be a finite subset of $\G$, and 
        $n$ words of three letters in  
        $\mathcal{A}$, $\tte_i= a(i)_1 a(i)_2 a(i)_3$,  with $i\leq n$,
        such that, for all $i$, 
        $\tte_i  \underset{\G}{=}   1$ in     $\Gamma$.

        Then, there exists an integer $l \in [1,
        10^{30}(\delta+1)^{10}]$, such that, 

     \[ \forall i \leq n, \;  \exists l_1,l_2,l_3 \in \tilG,   \;  \exists
     r_1,r_2,r_3 \in \tilG, \; \exists c_1,c_2,c_3 \in \tilG \]
             satisfying  
             \[  \forall j\leq 3,  \quad
             \mathfrak{rep}_l(a(i)_j) = l_j c_j r_j \quad r_j = (l_{j+1
               \,[3]})^{-1}. \]

             and either $\forall j\leq 3, \,  c_j$ has a long normal form of
             length at most $\kappa n$ with letters in  $ S \cup (\bigcup_{k=1}^q
             \Sec_k(\kappa))$, or  
             $\exists k, \forall j\leq 3, \, c_j \in P_k$.

        \end{theo}

Note that the elements $r_j$ and $l_j$ represent elements of the group $\tilG$, and therefore paths up to translation in $\widehat{Cay} \G$ (which makes a syntaxic difference with the situation in Theorem 3.1, where the cylinders were \emph{not} considered modulo translation).

  We say that a triangle $\tte_i$ is singular for its canonical
  representatives, if it falls in the second case of the theorem, and regular otherwise. 
  
  Let us describe the
  picture on the paths defined by $\frap$. Each regular equation gives
  a triangle $(1, z(i)_1, z(i)_1z(i)_2)$  that looks like a tripod, with a
  small defect at its center, that is neither long in perimeter, nor in maximal
  angle. Each singular equation gives a triangle $(1, z(i)_1, z(i)_1z(i)_2)$
  that is a tripod, whose center is a vertex fixed by a
  parabolic subgroup.

         \subsection{Canonical representatives are geometric.}

         We need geometric properties of canonical representatives in order to
         get some control on them.

         \begin{prop}
           Let $\gamma\in \G$, and $l$ an integer,  
           then its canonical representative
           $\mathfrak{rep}_l(\gamma)$ is geometric, 
           in the sense that it is in $\calL
           \subset \tilG$. 
         \end{prop}
         
         Recall (see the remark before Lemma \ref{lem;junctions}) 
         that no junction has an angle greater than $\theta$. Therefore if the
         path $\mathfrak{rep}_l(\gamma)$ has an angle greater than $\theta$,
         it must be between two junctions, that is on the center of a
         slice.  If it  pass again at this vertex with an angle greater than
         $\theta$, this means that two slices have same center, but this is
         impossible. Therefore $\mathfrak{rep}_l(\gamma)$ does not have any
         $\theta$-detour.

         Let us show that it is in fact a $(L_1,L_2)$-quasi-geodesic. Let $f$
         be an arc-length parameterization of the path.

     Since junctions are geodesics of length at most $1000  \delta$,
      it is enough to check the inequality 
      $|b'-a'| = length(f([a',b'])) 
      \leq 10^4\delta M  d(f(a'),f(b')) + (10^7\delta
      M -2000\delta)$, between centers of slices. 
      Let 
      $a'$ and $b'$
      be so that $f(a')$ and $f(b')$ are centers of slices in the slice
      decomposition. 
      Since any junction is shorter than  $1000\delta$,  
      the  path contains at least $ |b'-a'|/(1000\delta)$
      junctions  between these two points.  
      Therefore, $|b'-a'| \leq 1000\delta \times (n+1)$, 
      where $n$
      is the number of slices between $f(b')$ and $f(a')$.  

      Now, let us bound from below the distance between two points $x$ and $y$ 
      separated by $n$  slices in the decomposition of an arbitrary cylinder
      $Cyl(g,g')$. We need to refer to structural details in \cite{Dis} about
      cylinders, in order to compute this bound.
      First,  cylinders are union of cylinders that contain no parabolic
      slice, separated by very large angles (see 2.14 in \cite{Dis}). So it is enough to treat only 
      the case of cylinders without  
      parabolic slice, since their concatenation with large angle would remain quasi-geodesic.
 
      In these cylinders, by the definition of slices (see \cite{Dis} Def. 2.15), 
      $n\leq | \sharp (N_L^{g,g'} (x) \setminus
      N_L^{g,g'}(y) )   -    \sharp (N_L^{g,g'} (y) 
      \setminus N_L^{g,g'}(x) )  +
      \sharp (N_R^{g,g'} (y) \setminus N_R^{g,g'}(x) )   - 
      \sharp (N_R^{g,g'} (x)
      \setminus N_R^{g,g'}(y) )  | $, 
      where the set $N_L^{g,g'} (x)$ represents the set of 
      points of $Cyl(g,g')$
      at distance at least $100\delta$ from $x$ and that are closer to 
      $g$ than $x$ ($L$ stands for left, and $R$ for right: the role of $g$ and
      $g'$ are permuted). 
      Therefore, one of these four sets must have 
      cardinality at least $n/4$: say
      without loss of generality that $\sharp (N_L^{g,g'} (x) \setminus
      N_L^{g,g'}(y) )) \geq n/4$. 
      But a cylinder is contained in the union of cones 
      of radius and angle $50\delta$ centered at the edges of a 
      geodesic between $g$ and $g'$. 
      Therefore, as $M$ is an upper bound to the cardinality of cones of 
      radius and angle
      $50\delta$,  there are at least $n/(4M)$ different cones 
      centered on edges of a
      geodesic $[g,g']$ that intersect $(N_L^{g,g'} (x) \setminus
      N_L^{g,g'}(y))$. 
      Considering the extremal cones of this family, one gets two
      points  at least $n/(4 M)$ apart in this set. 
      One  concludes that $d(x,y) \geq n/(4M) - 200\delta$.

      With the minoration $1000\delta (n+1) \geq |b'-a'|$ previously  obtained, 
      one gets, for $x=f(a')$ and $y=f(b')$:  $d(f(a'),f(b')) \geq |b'-a'|/(4000\delta M) - 1/(4M)-  200\delta $ and thus    
      $|b'-a'| \leq 4000\delta M  
      \times  d(f(a'),f(b')) + 800000\delta^2 M +4000$. 
      $\square$

   \section{Lifting Solutions in $\calL$}
      
      We can now reproduce in the relative case the idea, presented in the
      hyperbolic case in section 1, that one can lift
      a solution of a system of equations as a solution of another system 
      in the group $\tilG$.

       Let $\mathscr{E}=\{\tte_1=1, \tte_2=1, \dots, \tte_n=1  \}$ be a finite
       system of triangular equations in $\G$, with set of unknowns $\Omega$,
       and set of parameters $P$. Let us write $\tte_i = z(i)_1 z(i)_2 z(i)_3$.

       Let us consider another instance of $\Omega$, which we distinguish:
       $\tilde{\Omega} \simeq \Omega$. Let $\Omega' = \tilde{\Omega} \cup 
       \{ l(i)_j, c(i)_j, r(i)_j \, |\,  i\leq n, j\leq 3 \}$. 
       
       Given $i\leq n$, 
       let us define the $6n$ equations $\widetilde{\tte}_{i,j}$, for $j\leq 6$:

       \[  \forall j\leq 3, \,  
       \widetilde{\tte}_{i,j}= (\widetilde{z(i)_j})^{-1}\, l(i)_j \,
       c(i)_j\, r(i)_j  \quad    \widetilde{\tte}_{i,j+3} =
       r(i)_j\,l(i)_{j+1\, [3]} \]

       We may or may not use the following equations, but we need to define them.
       
       \[ \forall i\leq n \; \widetilde{\tte}_{i,\hbox{c}} =
       c(i)_1\, c(i)_2\, c(i)_3 \]

       Let us now define the set $\mathscr{T}$ of triples 
     $(c_1,c_2,c_3)$ of $\tilG$, whose long normal form is of length at most
     $\kappa n$, contains only letters in $S \cup (\bigcup_{k=1}^q \Sec_k(\kappa))$
     and such that the concatenation of 
     $\frap(c_1), \frap(c_2)$ and $ \frap(c_3)$ is a closed loop, as in the
     statement of Theorem \ref{theo;repcan}. 
	In order to compute effectively this set  $\mathscr{T}$, it is enough, first, to compute   $\Sec_k(\kappa)$ for any $k$ and the explicit constant $\kappa$ (this can be done 
	knowing the list of simple loops of length $2\kappa$ in $\widehat{Cay}\G$), and second, to use a solution to the word problem in $\G$, which we have \cite{F}.

     We can now define a family $\mathfrak{F}$ 
     of systems of equations with parameters 
     in $\tilG$: a system of $\mathfrak{F}$ contains all  the equations
     $\widetilde{\tte}_{i,j} =1$, for $i\leq n, j\leq 6$, and, for all $i$, either
     it contains $\widetilde{\tte}_{i,\hbox{c}}$ (in which case we call
     the index $i$ singular), or it contains
     equations    $ c(i)_1 = c_1,  \, c(i)_2 = c_2, \, c(i)_3 =c_3$
     for some triple of parameters $(c_1,c_2,c_3) \in \mathscr{T}$ (in which
     case we call the index $i$ regular). With the preliminary knowledge of the items mentioned at the end of the introduction, since $\mathscr{T}$ is computable,   $\mathfrak{F}$ is computable.

     \begin{prop}\label{prop;lift}
       Let $\G$ be a torsion free relatively hyperbolic group, 
       and $\tilG$ as above.
       
       Given a finite triangular 
       system $\mathscr{E}$ of
       equations in $\G$, of set of unknowns $\Omega$,  there is a computable  
       finite family  $\mathfrak{F}$  of finite 
       systems of $\tilG$, of set of unknowns  containing a copy $\tilde{\Omega}$ of
       $\Omega$,  that satisfy the following.

       Given a solution $s$ of $\mathscr{E}$, there is a system $
       \tilde{\mathscr{E}}$ in  $\mathfrak{F}$  that admits a solution 
       $\tilde{s}$
       such that $\tilde{s}(\tilde{\omega})$ is in $\mathcal{L}$, and
       $\overline{\tilde{s}(\tilde{\omega})   } = s(\omega)$.

       If a system in $\mathfrak{F}$ has a solution $\tilde{s}$, then
       $\overline{\tilde{s}}$ restricted to $\tilde{\Omega}$ (pre-composed with
       the bijection between $\Omega$ and $\tilde{\Omega}$) is a solution to
       $\mathscr{E}$. 
     \end{prop}

     We choose  $\mathfrak{F}$ as before the statement of the proposition.  
     As already said, $\mathfrak{F}$ is computable.

     In order to prove the first assertion, 
     let us first choose correctly $\tilde{\mathscr{E}}$ (note that it may not be computable a priori, however we only need its existence). 
	We need to choose
     which indices are singular, and which are regular. Given $s$, we consider
     the canonical representatives of the solution given by Theorem
     \ref{theo;repcan}. We say that an index $i$ is singular if and only if
     the corresponding triangle is singular, as defined in the theorem.

     With the notations as in the beginning of the section, let $\omega=
     z(i)_j$. We then assign  the values for $l(i)_j, c(i)_j, r(i)_j$ for 
     $ i\leq n, j\leq 3$ as given by Theorem
     \ref{theo;repcan} for the element $\mathfrak{rep}_l(z(i)_j)$. 
     This theorem precisely tells that this is a solution
     to the corresponding system of equations in $\mathfrak{F}$.

     To check the last assertion, 
     simply consider the image in $\G$ of the solution in $\tilG$, by
     definition, this is a solution of $\mathscr{E}$. $\square$

   \section{Results, and Conclusion}

      We make use of important results about decidability of existential theory 
      for free products of groups. We are interested in results of V.~Diekert, M.~Lohrey, and A.~Muscholl. 
 In the case our parabolic groups are free abelian, we considered a free product of abelian subgroups, and the following is sufficient for our study (recall that in our context, generalized normalized regular languages are normalized regular).

      \begin{theo}(Diekert, Muscholl, \cite[Thm 25]{DiMu}) \label{theo;DM}

        The existential theory with normalized rational constraints of 
        a free product of abelian groups is decidable.
      \end{theo}

      For general free products, we use a special case of Theorem 3.10 in \cite{DiLo} (or Theorem 4 in \cite{DiLo_a}), for the class of generalized rational constraints that we introduced above:

      \begin{theo}(Diekert, Lohrey, \cite{DiLo}) \label{theo;DL}

        For
        any free product of groups whose existential theory 
	with parameters is decidable, 
        the existential theory, with generalized normalized rational 
	constraints (with respect to the given free product),
        is decidable.
      \end{theo}

     We obtain the following  result 
     which is     slightly stronger than needed for proving that existential
     theories  are decidable. However, while discussing with D.~Groves
     on how to apply these results to the Isomorphism Problem, it became
     apparent that it would be useful for this purpose at least. 
     It is the analogue of Proposition \ref{prop;isom_hyp}
     that was stated in the hyperbolic case.

     \begin{prop} \label{prop;isom}

              Let  $\G$ a  torsion free relatively hyperbolic group, and $\tilG$ the
              free product of its parabolic subgroups and of a free group, as
              above. Let $\calL$ the set of geometric elements of  $\tilG$ defined before.

              There is an algorithm that performs the following.

              Given a finite system $\mathscr{E}$ of triangular equations in $\G$, of
              set of unknowns $\Omega$,  and given,
              for each  $\omega \in \Omega$, a $\mathfrak{gnr}$  language  
              $\mathcal{L}_{\omega}$  of geometric elements ({\it i.e. } 
              such that $\mathcal{L}_{\omega} \subset \calL$), the algorithm
              always terminates, and says ``yes'' if there exists a solution
              $s$
              of $\mathscr{E}$ with the property that, for all $\omega$,  
              \emph{any}
              pre-image of  $s(\omega)$ in $\calL$ is in
              $\mathcal{L}_{\omega}$, and only if  there exists a solution
              $s$ of $\mathscr{E}$ with the property that, for all $\omega$,  
              \emph{some}
              pre-image of  $s(\omega)$ in $\calL$ is in
              $\mathcal{L}_{\omega}$

              Moreover if  $\G$ has only free abelian parabolic subgroups, the same is true with 
              $\mathcal{L}_{\omega}$ normalized regular (not necessarily generalized).
            
     \end{prop}

     {\it Proof. } Let $\tilde{\Omega}$ be a set in bijection with $\Omega$,
      that we interpret as a set of unknowns in $\tilG$. 
      Let us define the set of constraints $\widetilde{\mathscr{C}} = \{
      \tilde{\omega} \in \calL_{\omega}, \, \forall  \tilde{\omega} \in \widetilde{\Omega}
      \}$, where $\omega$ is the element of $\Omega$ in bijection with
      $\tilde{\omega}$.

        Then we apply Proposition \ref{prop;lift}, to  compute a finite  family 
        $\mathfrak{F}$ of systems
      of equations in $\tilG$ that satisfy then the following (namely the first
        and second conclusions of Proposition \ref{prop;lift}).
 
      First, if $s$ is a solution of $\mathscr{E}$  with the property that, for all $\omega$,  
              \emph{any}
              pre-image of  $s(\omega)$ in $\calL$ is in
              $\mathcal{L}_{\omega}$, then there is a system of $\mathfrak{F}$
              with a solution $\tilde{s}$ such that
              $\tilde{s}(\tilde{\omega})$ is a pre-image of $s(\omega)$ in
              $\calL$, hence by assumption, in $\mathcal{L}_{\omega}$. 

      Second, if $\tilde{s}$ is a solution of some system of $\mathfrak{F}$,
      with constraints  $\widetilde{\mathscr{C}}$, then $\overline{\tilde{s}}$
      induces a solution of  $\mathscr{E}$, and clearly, for all $\omega$,
      some pre-image of $\overline{\tilde{s}}(\omega)$ is in $\calL_{\omega}$ 
      (namely $\tilde{s}(\omega)$).

     Thus, it is enough to have an algorithm that solves the following
     question: is there a system among those in  $\mathfrak{F}$ that has a
     solution in $\tilG$, with
     the constraints $\widetilde{\mathscr{C}}$ ? 
     
     This is exactly what Theorem \ref{theo;DL} provides in the case of
     $\mathfrak{gnr}$ languages.
     
      In the case of free abelian parabolic subgroups and normalized regular
      languages, 
      $\tilG$ is a free product of abelian groups, 
     and  Theorem \ref{theo;DM} can be applied instead to give the result.   $\square$

      \paragraph{Proof of Theorem \ref{theo;intro2}.}

      Let $\mathscr{E}$ be a finite system of equations, and $\mathscr{I}$ be
      a finite system of inequations. As already discussed in the first
      section, 
      one can assume that all the equations are triangular, and that  $\mathscr{I} = \{ \omega
      \neq 1,  \omega \in \Omega_i\}$, for some $\Omega_i \subset \Omega$.

      We then apply Proposition \ref{prop;isom} for $\calL_\omega= \calL$ when
      $\omega \notin \Omega_i$, and  $\calL_\omega=\calL\setminus \calL_0$ when
      $\omega \notin \Omega_i$, which are generalized 
      normalized regular languages (with explicit  $\mathfrak{gnr}^\bullet$-automata)  
      by  Proposition \ref{prop;reg}, 
      and Corollary \ref{coro;L_0}. 

      Let us now notice that a possible solution $s$ has one of its lifts in
      $\calL$ satisfying these constraints if and only if it has all of them,
      since a lift of $s(\omega)$ belongs to $\calL_0$ if and only if
      $s(\omega) = 1$.

      Therefore, the algorithm provided by  Proposition \ref{prop;isom}
      answers the question whether  $(\mathscr{E}, \mathscr{I})$ has a
      solution. $\square$

      \paragraph{Proof of Corollary \ref{coro;intro3}.}

      The corollary is implied by Theorem \ref{theo;intro2} and the fact that
      existential theory of finitely generated virtually abelian groups 
      is decidable. This latter
      property is folklore, but note that it is not immediate 
      (for instance it
      is unknown whether virtually free groups have decidable 
      existential theory, although, needless to say at this point,  
      it is known for free groups). 
      Thus we give a proof here.

      Any such group $\Gamma$   
      is abelian-by-finite: it splits as
      an exact sequence $1 \to A \to \Gamma \to F \to 1$ 
      where $A$ is free abelian of finite rank,
      and $F$ is finite. Hence we are left with the lemma.

      \begin{lemma}
        There exists an algorithm satisfying:

        Input:  an integer $n$, a finite group $F$, a morphism 
        $\rho \colon F \to GL_n(\mathbb{Z})$,  a finite presentation of 
        $\Gamma$   splitting  as $1 \to \mathbb{Z}^n \to \Gamma \to F \to 1$,
        with the given action of $F$ on $\mathbb{Z}^n$,  
        a (set theoretic) cross section $s: F\to \Gamma, s(F) = \bar{F}$, 
        and a system of equations and inequations in  $\Gamma$. 

        Output:  whether there exists a solution of the system in $\Gamma$.  
      \end{lemma}

      Given a finite system of equations and inequations in $\Gamma$, one can
      project the system in $F$ and find all the solutions. For each solution
      $(f)$ of the projected system, one can lift $(f)$ by the cross section
      $s$  in $\Gamma$ in the preferred set $\bar{F}$ of representatives of
      $F$, in order to get a system in $\Gamma$ where unknowns are in
      $A=\mathbb{Z}^n$.  For example, if an equation of the original system is
      $\phi_j: z_{(j,1)}\, z_{(j,2)}\, z_{(j,3)} = 1$, one changes it into 

      $$
      \phi'_j:  \left( \bar{f}_{(j,1)} a_{(j,1)} \right ) \,  
      \left( \bar{f}_{(j,2)}a_{(j,2)}\right) \,
      \left( \bar{f}_{(j,3)}a_{(j,3)}\right)\: = \:1,
      $$ 
      where $\bar{f}_{(j,i)} \in \bar{F}$ are the images by 
      $s$ of the solutions in $(f)$, (hence are known
      parameters),   and where  $a_{(j,i)}$ is  a parameter 
      if $z_{(j,i)}$ is a
      parameter,  or unknown, if $z_{(j,i)}$ is an unknown,  in $A$.
      One gathers
      all the $\bar{f}_{(j,i)}$ to the left, by using the action of 
      $\rho(F)<GL_n(\mathbb{Z})$ on $A$: for all $f  \in \bar{F}$,
      $ a f=f \,  (\rho(f) (a))$ for all $a
      \in A$.  
      As we know explicitly the map $\rho$, one can compute 
      $a'_{(j,i)}$ satisfying
      $\bar{f}_{(j,i+1)}\, a'_{(j,i)} =a_{(j,i)} \, \bar{f}_{(j,i+1)} $,  
      in terms of the
      coefficients of $a_{(j,i)}$ expressed in our basis of $A$.  
      Therefore one is
      left with the equation 
      $$
      \left( \bar{f}_{(j,1)}\bar{f}_{(j,2)}\bar{f}_{(j,3)}\right) \,
      \left( a''_{(j,1)}a'_{(j,2)}a_{(j,3)} \right)  \: =\: 1
      $$ 
      (and similarly for inequations). 
      The first product of three terms is assumed to be
      in $A$, since we started with a solution of the projected system in $F$. 
      The results of these products of three elements in 
      $\bar{F}$ that are in $A$ (there are finitely many) 
      are computable, thus one is left with a 
      system of linear equations in the 
      coefficients of the $a_{(j,i)}$ in a chosen basis of $A$.  
      Such systems are decidable, and, by construction, 
      the original system has a solution if, 
      and only if, one of the new linear 
      systems (associated to solutions of the projected system in $F$) 
      has a solution.   $\square$

%   FIN DE LA REVISION ... VIRER LE RESTE... GL   ET HRUSKA

\part{Finding relative hyperbolic structures}

The aim of this note is to provide a  general algorithm that recognizes a class of
groups, and that allows to  unify   
a certain number of algorithms
originally designed for every group in this class. 

The class in question is that of relatively hyperbolic groups. It was introduced by 
M.~Gromov \cite{Grom}, and there is now a rich growing literature on 
this subject (see works of B.~Bowditch, C.~Dru\c{t}u, A.~Yaman, D.~Osin, D.~Groves,
myself, and many others).  
It generalizes the class of geometrically
finite Kleinian groups, in Gromov's hyperbolicity spirit.

In order to run a certain number of known algorithms in these groups ({\it
  e.g.} B.~Farb's solution to word problem \cite{F},
I.~Bumagin's solution to conjugacy problem \cite{Bu}, 
D.~Osin's solution to the root problem \cite{O}, my solutions to Diophantine and
existential theory decidability presented in part 1 of the present e-print, etc.),
one needs to first find correct values for some characteristic constants, such
as hyperbolicity constant. Our algorithm will provide a way to do this, when
given  a candidate presentation compatible with the relative structure. 
 Better: we can  detect 
whether a given finitely presented group is hyperbolic relative to 
abelian subgroups,  and find the structure in such case. Thus, for each
  algorithmic problem mentioned, a single algorithm solves it for all the
  groups  hyperbolic relative to 
abelian subgroups (and torsion free for those in part 1).

  This was inspired by the
algorithm of P.~Papasoglu for recognizing hyperbolic groups, and
motivated by  the author's collaboration with D.~Groves to study the isomorphy problem (also known as isomorphism problem) for
 relatively hyperbolic groups, in which our algorithm is an important tool.

\begin{theo}\label{theo;intro1bis}
There is an algorithm such that, given a finite presentation of a
group $\Gamma$, the generators in $\Gamma$ of finitely many subgroups
$H_1,\dots, H_n$, and a solution to the word problem in each of the $H_i$ with
the given generators, terminates if and only if
 the group is hyperbolic relative to $H_1, \dots, H_n$.

In this case, the algorithm provides an explicit bound for the linear
relative isoperimetric inequality satisfied by the relative presentation. 

\end{theo}

We give explicitly such an algorithm, and we use it in Corollary \ref{coro;delta}, to find explicitly
the hyperbolicity constant and the list of simple loops of given length of
Farb's coned-off graph.

In the case of  abelian parabolic subgroups, we propose here a sharper
result. 
  Again, the algorithm  is explicitly given.

\begin{theo}\label{theo;intro2bis}
 There is an algorithm whose input is a finite presentation of a group,  that terminates 
  if and only if 
  the group is hyperbolic 
       relative to abelian subgroups, and gives as  output  
       an exact relative presentation. 

In particular, the languages of the presentations of  relatively hyperbolic
groups with  abelian parabolics, and of those with free abelian parabolics are recursively enumerable.

\end{theo}

This latter property features in Sela's list of problems for  limit groups
(which I proved hyperbolic relatively to free abelian parabolics in previous work).

\section{Relative presentations, and hyperbolicity}
\subsection{Lengths and areas} 

We present here the framework, which is inspired by \cite{O} and \cite{Rebb}.

For any group $H$, we denote by $\calT(H)$ its table of multiplication 
of $H$,  $\calT(H) = \{ ((a,b,c) \in H^3, abc =_{H} 1   \}$. 

We denote by $\Gamma$ a finitely generated group, and $X$ a finite symmetric set of
generators. Assume that we are given $n$ subgroups $H_1, \dots, H_n$ of
$\Gamma$, and,  for all $i=1, \dots, n$, 
finitely generated groups $\tilHi$, with epimorphisms $\pi_i:\, \tilHi \tto H_i$.    We define the free product $F= F_X * \tilde{H}_1 *  \dots * \tilde{H}_n  $, where $F_X$
is the free group on $X$.  We define  $\calX = X \sqcup \tilde{H}_1 \sqcup  \dots \sqcup \tilde{H}_n$, and 
 $\calS = \bigsqcup_{i=1}^n \calT(\tilHi)$. By abuse of language, we do not make any distinction between a triple
 $(a,b,c)$ in $\calS$, and the word $abc$ in the alphabet $\calX$.

We say that the group $\Gamma$ has a finite $(\tilHi)_{i\leq n}$-relative  presentation   
 if there exists a finite set $\calR \subset \calX^*$ of words over $\calX$ (which can, without loss of generality , be assumed to consists of words of $2$ or $3$ letters)
such that  $ F/\langle\langle \calR \rangle \rangle$  is isomorphic to    $\Gamma$, by a morphism that is identity
on $X$ and induces $\pi_i$ on each of the $\tilHi$.

One easily checks that, in this case, $ \Gamma \simeq  \langle \calX \; | \; \calR \cup \calS \rangle$ (thus $\calX$ is seen as generating set of $\Gamma$). The elements of $\Gamma$ can be measured with the word metric on $X$ (which we  call the absolute length) and also with the word metric on $\calX$ (which is the relative length).

Let us notice  the finite subset $\calX_\calR$ of $\calX$ consisting of the  elements of $ X$ and those in some  $\tilHi$ appearing (as letters) in some word in  $\calR$. They have relative length $1$ but their absolute length can be large.

We say that the finite relative presentation is \emph{exact} if  $\tilHi \simeq H_i$ for all $i$, and each $\pi_i$ is the identity.

{\it Remark. } In general we say that $\Gamma$ is \emph{finitely presented
  relative to $H_1, \dots, H_n$} if it admits an exact finite relative presentation (see for instance \cite{O}). If a group admits an non-exact finite relative presentation, then it admits an exact one, and if  $\Gamma$ is finitely presented it has a finite relative presentation  for every family of groups $\tilHi$ that maps onto $H_i$.

Let us temporarily consider more generality. Given a group presentation
$G=\langle S\, |\, R \rangle$, a Van-Kampen diagram for a relation $w\in S^ *$, $w=_G 1$,  
in the  alphabet $S$,  is a finite, planar, simply connected,  2-complex whose
 oriented edges are labeled by elements of $S$, such that the word labeling
 of the boundary of any
 cell is (up to a cyclic permutation, and choice of orientation) 
 an element of $R$, and such that the boundary of the complex consists of one
 loop, labeled by the word $w$ (up to cyclic permutation).

The \emph{combinatorial area} of a relation $w\in S^*$ of $G$ ({\it i.e. }
$w=_G 1$) is the minimum, over the set of Van Kampen diagrams of $w$ for the
presentation $G  =   \langle S \,|\, R  \rangle $, of the number of cells of the diagram.

In the following, we use these notions for the relative (not necessarily exact) presentation given.

 The following definition is coherent with other definitions of relative hyperbolicity.

\begin{defprop}(Osin \cite{O}, Rebbechi \cite{Rebb}) 

  A group $\Gamma$ is hyperbolic relative to a family of subgroups $H_1, \dots, H_n$ if, and only if, 
  it admits an  exact relative presentation that satisfies a linear isoperimetric inequality ({\it i. e } there exists $K>0$, such that  for all relation $w$ of  $\Gamma$ the combinatorial area of $w$ for the relative presentation is at most $K$ times the relative length of $w$). 

  In this case, we call the subgroups $H_1,\dots,H_n$ and their conjugates, the parabolic subgroups of $\Gamma$.

\end{defprop}

\subsection{Clusters in Van Kampen Diagrams}

  In this paragraph, $\Gamma$ is a group with a (non necessarily exact) relative presentation, and 
  $w \in \calX^*$ is a relation of $\Gamma$ ($w =_\Gamma 1$). The lengths and combinatorial areas are measured in this relative presentation.
We intend to gather some useful results on the structure of Van Kampen diagrams for $w$, with cells 
  labeled by elements of
  $\calR \cup \calS$.

  Given $i\leq n$, a $i$-\emph{cluster} of a Van Kampen diagram $D$ is a sub-complex of $D$ that contains
  only cells labeled by elements of $\calT(\tilHi)$, that is connected, 
  without
  global cut point,    
  and that is maximal for this property. Sometimes we do not specify the index $i$.
  For a cluster $C$, the boundary $\dron C$ of $C$ consists of its 
  edges that belong to only one cell of $C$.  The \emph{exterior boundary} of
  $C$ is the unique component of $\dron C$ that is homotopic to $\dron D$ in
  $D\setminus C$.

  \begin{lemma}\label{lem;inner_edges}(Inner edges of $D$)

    Any edge of a Van Kampen diagram belongs to at most one of its clusters.

    For any cluster $C$ of a Van Kampen diagram $D$, any edge in $\dron C$ 
    is either in $\dron D$, or is in another cell of $D$ labeled by an element of $\calR$.  
  \end{lemma}
  
  \begin{proof} 
If an edge $e$ belongs to a $i$-cluster $C$, then the label $s$
  of $e$ belongs to $\tilHi$. If $C'$ is a  $j$-cluster
  containing $e$, $s\in \tilde{H}_j$, and since in $\calX$, 
        $\tilde{H}_i$ and $\tilde{H}_j$ are disjoint or equal, 
   one has  $j=i$, and by maximality, $C=C'$.
  
  Thus an edge of $\dron C$ does not belong to any other cluster. If it  is
  not in $\dron D$, it  belongs
  to another cell in $D$, that must be therefore labeled by an element of $\calR$. \end{proof}

\begin{lemma}\label{lem;find_minimal}(Small clusters in a non minimal diagram)

  Let $A>0$, and $w$ be a relation of $\Gamma$.  
        Assume that $w$ labels the boundary of a (non necessarily minimal) 
        Van Kampen diagram 
        $D$ with less than $A$ cells. Assume furthermore that each 
        cluster 
 contains at most one edge in the boundary $\dron D$.

        Then,   $w$ is labeled by letters of $\calX_\calR$,  or letters in some $\tilHi$ that are product, as elements of this group,  of at most $3A$ elements of $\calX_\calR \cap \tilHi$. 
   Moreover, there exists a minimal Van Kampen diagram for $w$ with 
        at most  $A$ cells,  
        whose every cluster has
        perimeter at most $3A$.

\end{lemma}
 
 \begin{proof}  Consider a letter $s$ in $w$, and  consider the edge $e$ it labels in $\dron D$.  Either it belongs to a cell of $D$ labeled by an element of $\calR$,  or it belongs to some cluster, and then $e$ is homotopic in $D$ to a path $b$ such that the concatenation $(b e)$ is a boundary component of the cluster. Now $b$ is assumed to contain no edge in $\dron D$, thus, according to Lemma \ref{lem;inner_edges}, only edges that are labeled by elements of $\calX_\calR$. Moreover, since there are at most $3A$ different edges in $D$, the length of $b$ is at most $3A$. This shows the first assertion.

Now consider a minimal Van Kampen diagram $D_m$ for $w$. It contains at most
$A$ cells, by assumption. Any cell not labeled by an element of $\calR$ is in
a cluster, by definition. 
       The length of  boundary components of a cluster cannot exceed
      the total number of edges of $D_m$, that is at most $3A$. 
      This proves the lemma.  \end{proof}

In order to prove the next proposition we will use several times the following
observation.

\begin{lemma}\label{lem;Euler}
  Let $D$ be a  planar simplicial 2-complex, and let $b$ be a simple loop of
  length $l(b)$, in its
  1-skeleton. If $D_b$ is a simply connected sub-complex, 
  such that any edge of $b$ is contained
  in a triangle of $D_b$, then its area is at least $l(b)-2$.
\end{lemma}
 
   \begin{proof}   The Euler characteristic of $D_b$ is $\chi(D_b) = 1$  and  
  if $D_b$ contains $N$ triangles,
  it contains $3N/2 + l(b) /2$ edges. If the number of vertices that are not in
  $b$ is $x$, one gets $\chi(D_b) = N - (3N/2 + l(b) /2) + l(b) +x$. Hence $N=l(b)-2 + 2x \geq l(b)-2$. \end{proof}

\begin{prop}\label{prop;reduction} (Simplifying examples of high ratio (area/perimeter))

Let $A>0$. 
  Assume that,  for all $i$,   all $r\leq 3A$, and all  $s_j \in \calX_\calR \cap \tilHi$, for $j=1, \dots, r$, 
the word  $s_1 \dots s_r$ is a relation of $\Gamma$ only if it is a relation of $\tilHi$.

  For any relation $w$, with $Area(w)\leq A$, and $\frac{Area(w)}{|w\absx }>3$, 
  there exists a 
  relation $w'$ such that:

\begin{itemize}

\item  one has $\frac{Area(w')}{|w'\absx} \geq \frac{Area(w)}{3|w\absx }, $

\item every letter of $w'$ is either in $\calX_\calR$ or is a letter in some $\tilHi$ that is product, as element of this group,  of at most $3A$ elements of $\calX_\calR \cap \tilHi$.

\item there exists a minimal Van Kampen diagram for $w'$ with at most  $A$ cells, each of which  either is labeled by an element of $\calR$, or is in a simply connected cluster of perimeter at most  $(A+2)$.

\end{itemize}

\end{prop}

  \begin{proof}   Let $w$ be a relation, and $D$ a Van Kampen diagram for $w$,
that is minimal for the area, and among all possible minimal diagrams, that is
minimal for the number of cells labeled by $\calR$. Let $C_1, \dots, C_r$ be the
clusters of $D$ that contain two edges of $\dron D$. Clearly $r\leq |w\absx
/2$.   Let $\{a_1, \dots, a_s\}$ the collection of all the maximal arcs in the
boundary of these clusters, joining $\dron D$ to $\dron D$,  
that do not contain an edge in $\dron D$. Each arc
$a_j$ is labeled by a certain word $h_j$ with letters in some $\tilHi$. Let
$s_j$ be the element of $\tilHi$ defined by $h_j$.

An arc $a_j$ cannot cross a cluster, it follows easily that there are at most
$|w\absx-1$ such arcs.  A \emph{piece} of $D$ is the closure of a 
component of $D \setminus \bigcup_{j=1}^s a_j$.
Thus a piece is a connected sub-diagram of $D$, and there are thus at most $s+1\leq |w\absx $ different pieces.

Let $p$ be a piece of $D$, and $b$ a boundary component of $p$. 
It is labeled by $w_p= w_1 h_{i_1} \dots w_{n(p)} h_{i_{n(p)}}$ where the $w_i$ are sub-words of $w$,
corresponding to arcs of $\dron D \cap b$,  and the $h_i$ are as previously defined.

We then define the word $w'_p= w_1s_1\dots w_{n(p)} s_{n(p)}$. Note that
$|w'_p\absx$ is the sum of length of sub-words of $w$ involved, plus
$n(p)$. Moreover, since $D$ is planar, each arc $a_i$ appears in exactly two
pieces, and therefore the sum of the perimeter of all the $w'_p$  is  $|w\absx + 2s \leq 3 |w\absx$.

 Let $D_p$ a minimal Van Kampen diagram for $w'_p$.  Gluing all the
diagrams $D_p$ along the edges corresponding to the letters $s_j$ provides a
Van Kampen diagram for $w$, and therefore $Area(w) \leq \sum_{p} Area (w'_p)$.

If for all $p$, $Area(w'_p) <  \frac{Area(w)}{3|w\absx }  |w'_p\absx $,
 then $Area (w) <  \frac{ Area(w)}{3|w\absx } \sum_{p}
 |w'_p\absx$, which contradicts that the sum of perimeter of pieces is at
 most $3 |w\absx$. Therefore, there is a piece such that 
$Area(w'_p) \geq \frac{Area(w)}{3|w\absx} |w'_p\absx$, which is   $ >
 |w'_p\absx $ by assumption.

Assume that $p$ is the component of a cluster $C$.  We will show that $C$ is a
  simply connected cluster, and then derive a contradiction. Assuming the contrary, 
   there is a boundary component $b$ of $C$, different from the exterior
  boundary of $C$, 
  labeled by letters of
  $\calX_\calR$,  and its length is at most the total number of edges, that is
  at most $3Area(w)\leq 3A$.   But by
  assumption, a word $w_b$ labeling such a boundary $b$ is a relation of the group
  $\tilHi$. One gets easily a Van Kampen diagram for $w_b$ with $|w_b\absx-2 (=l(b)-2)$
  cells labeled by elements of $\calT(\tilHi)$. 
  On the other hand,
  if $D_b$ is the minimal simply connected sub-diagram of $D$ bounded
  by $b$,  by Lemma \ref{lem;Euler} is has area greater than $l(b)-2$.
  Therefore, one can replace, without
  changing the area of $D$, and decreasing the number of cells labeled by
  $\calR$, the sub-complex $D_b$ by another whose cells are labeled by
  elements of $\calT(\tilHi)$. Since $D$ was chosen minimal for these
  properties, this means that all cells of $D_b$ belong to a cluster, which by
  Lemma \ref{lem;inner_edges} is necessarily $C$, 
  and this contradicts the maximality of the cluster $C$. 
  Thus, the cluster $C$ is simply connected, and therefore equals $p$. 
  But then, its boundary is labeled by a relation of $\tilHi$, and by
  minimality of the area, the number of  cells of the piece $p$ is at most its perimeter, contradicting the assumption on $p$.
Therefore, $p$ is not the component of a cluster.

Let us now propose a (not necessarily minimal) Van Kampen diagram for
$w'_p$. Start with the diagram $p$ itself. By assumption on the maximality of
the family $C_1, \dots, C_r$, there is no cluster of $p$ that has two edges on
$\dron p$. Moreover, for all $i$, such that $a_i$ is an arc of $\dron p$, each
edge of $a_i$ belongs to a cluster of $D$ that has no cell in $p$. Therefore,
no cluster of $p$ has an edge in any of the arcs $a_i$ in $\dron p$.   
   Now, for each $i$  for which  $a_i$ is an arc of $\dron p$,  by
   definition of $s_i$ the word $h_i(s_i)^{-1}$ is a  relation of $\tilHi$. 
   Hence, there exists a diagram $\delta_i$ bounded by this word, consisting
   of  at most $|h_i\absx$ cells, all of them  labeled by elements of some
   $\calT(\tilHi)$.     Let us define  $D_p$ to be the diagram obtained by
   gluing each $\delta_i$ on $p$ along the arc $a_i$. This is a diagram of
   $w'_p$, and its clusters are exactly the $\delta_i$, and the clusters of
   $p$. In any case, each of them contains only one edge in $\dron D_p$.

   Now $Area(D_p) \leq Area(p) + \sum_i Area(\delta_i)\leq Area(p) + \sum_i
   (|h_i\absx -1)$. 
 The exterior boundary $b_i$ of the cluster $C_i$  is a simple loop 
 longer than $|h_i\absx +1$, thus by Lemma \ref{lem;Euler}, 
 if $C'_i$ is a minimal
 simply connected diagram bounded by $b_i$, $Area(C'_i) \geq  |h_i\absx -1$.
  Therefore $Area(D_p) \leq Area(p) + \sum_i Area(C'_i)$, and one gets that
   $Area(D_p) \leq Area(D)\leq A$.

   Therefore, the relation $w'_p$ satisfies the assumption of the Lemma
   \ref{lem;find_minimal}, hence its letters have the requested form, and  there exists a minimal Van Kampen diagram
   $D'_p$ for $w'_p$ with at most $A$ cells either labeled by elements of
   $\calR$,  or in a cluster 
   whose boundary components are labeled by at most 
   $3A$ letters.

 It remains to see that every cluster can be assumed simply connected, and of perimeter $\leq A+2$. For
  that we proceed as before. 
Given a cluster $C$ of $D'_p$, denote by $b$ its exterior boundary, and by
  $D'_b$ the minimal simply connected sub-complex bounded by $b$. 
  By Lemma \ref{lem;Euler}, it contains at least $l(b)-2$ triangles.
  On the other hand,  since $l(b) \leq 3A$, by assumption of the lemma, the loop $b$ can be filled with  only cells labeled by elements in some $\calT(\tilHi)$, and clearly, only $l(b)-2$ such cells are sufficient. Therefore,  one can assume, without enlarging the area of $D'_p$, that every  cluster is simply connected. By Lemma \ref{lem;Euler}, its perimeter is at most $A+2$.
  This proves the lemma.  \end{proof}

\section{Algorithms}
\subsection{Algorithm of recognition.}

      In \cite{Papa}, P.~Papasoglu proves the following result, that appears
      in a different form in B. Bowditch's work (\cite{Bvp} Proposition 8.7.1)
      (see also M. Gromov \cite{Grom}). We choose to state it in the 
      context of group theory, rather than hyperbolic spaces, so we use the
      formulation of Papasoglu.

      \begin{prop}(\cite[Main theorem]{Papa})\label{prop;papa}
        Let $G$ be a group, with a triangular 
        presentation $G=\langle S|R\rangle$. Let $K >0$. 
        If there exists a word $w$ in $S$ that represents the trivial 
        element in 
        $G$, with combinatorial area $A(w)>K l(w)$, 
        then there exists another word   $w'$ with 
        combinatorial area $A(w') >\frac{1}{2.10^4} l(w')^2$ and 
        $A(w') \in [\frac{K}{2},240K]$.  

      \end{prop} 

      The proof takes place in a single Van Kampen diagram, with triangular
      cells, and thus does not need finiteness of the presentation.

\begin{lemma}\label{lem;reduction_finale}
      Let $K>0$.
 Assume that,  for all $i$,   all $r\leq 3(240K)$, and all  $s_j \in \calX_\calR \cap \tilHi$, for $j=1, \dots, r$, 
the word  $s_1 \dots s_r$ is a relation of $\Gamma$, only if it is a relation of $\tilHi$.

   Then,        either the presentation satisfies a linear isoperimetric
      inequality of factor $K$, or there is a
      relation $w''$ of area at most $240K$ such that 

\begin{itemize}
 \item one has       $\frac{Area(w'')}{|w''\absx} > \frac{\sqrt K}{600}$  

 \item every letter of $w''$ is either in $\calX_\calR$ or is a letter in some $\tilHi$ that is product, as element of this group,  of at most $3\times 240K$ elements of $\calX_\calR \cap \tilHi$.

\item and $w''$ admits a minimal Van Kampen diagram 
        whose every cell is either labeled by an element of $\calR$, or is contained in a simply connected cluster of perimeter at most 
        $(240K +2)$.
\end{itemize}
 
In particular, the letters labelling the boundaries of the clusters of a minimal Van Kampen diagram of $w''$   are equal to products of at most $3\times 240K$ elements of some $\calX_\calR \cap \tilHi$.

\end{lemma}
      
  \begin{proof}  If there is a relation $w$ of area
       $Area(w)>K |w\absx$, then by Proposition \ref{prop;papa}, 
       there is a relation $w'$ with area in
       $[\frac{K}{2},240K]$, and such that $Area(w') >\frac{1}{2.10^4}
      |w'\absx^2$. Therefore, by Proposition \ref{prop;reduction}, there is a
       relation $w''$ of area at most $240K$ such that
       $\frac{Area(w'')}{|w''\absx} \geq (Area(w')/|w'\absx^2)^{1/2} \times
      Area(w')^{1/2} /3$, and whose letters have the suitable form, and that 
	has a suitable Van Kampen diagram.  
      One gets that $\frac{Area(w'')}{|w''\absx} >  \sqrt{ 1/(2 . 10^4)} \times \frac{\sqrt{K/2}}{3} \geq   \frac{\sqrt K}{600}$. 

The last assertion is then a consequence of Lemma \ref{lem;inner_edges}.  \end{proof}

\begin{prop}\label{prop;algo}
There is an algorithm that has 
  the following properties.

The input is a non necessarily exact relative presentation, of a finitely generated group $\Gamma$,   relative to some groups  $\tilde{H}_1, \dots, \tilde{H}_n$: specifically, the finite generating set $X$ of $\Gamma$ is given, generators of each $\tilHi$ are given, the sets $\calX_\calR$, and $\calR$ are given, solutions to the word problem of each $\tilHi$ (for the respective generating sets)  are given.  Moreover, as input the algorithm takes either a solution to the word problem of $\Gamma$, or a blank entry. 

Assume that, either the presentation given is exact, or a solution to the word problem of $\Gamma$ has been given.

\begin{itemize}

\item Either the algorithm detects that the presentation given is non exact, or it stops if and only if it  satisfies an isoperimetric
inequality, and provides  a correct factor  $K$ such that, for all relator $w$,
$Area(w)\leq K |w \absx$.

\item If the presentation is non-exact, and such that one quotient $\tilHi
  \tto H_i$ has infinite kernel, and if a solution to the word problem for
  $\Gamma$ is given,  then the algorithm stops and tells that the presentation is non-exact.
\end{itemize}

\end{prop}

Note that Theorem \ref{theo;intro1bis} is a direct consequence of the Proposition.

  \begin{proof}    Our algorithm is as follows.
      Let $K$ be an integer variable.
        While the algorithm is not stopped, do the following instructions in order.

        \begin{enumerate}

        \item   If no solution of the word problem of $\Gamma$ is given in
        input, skip this instruction. Otherwise,  make, for all $i=1,\dots,
        n$,  the list of all words of at most $3(240K)$ letters in
        $\calX_\calR \cap \tilHi$, that are relations of $\Gamma$ (the
        solution of the word problem for $\Gamma$ provided in input allows to
        compute this finite list). With the solution to the word problem for
        $\tilHi$ given as input, check whether  these words are relations of
        $\tilHi$. If all of them are, proceed to instruction 2. If one of them
        is not, stop the algorithm, and return ``The presentation is not exact''.

        \item  For all $i=1,\dots, n$, make the list $\calV_i \subset \calX$ of all products in $\tilHi$  of at most $3(240K)$ elements of  $\calX_\calR \cap \tilHi$. This list is finite and easily computable. Make the list  $\calV_i'$ of all the words  of at most $(240K+2)$ letters in  $\calV_i$ that are trivial in the group $\tilHi$ (this requires a solution to the word problem in $\tilHi$, which is provided).

        \item  For all $i=1, \dots, n$, make the finite list $\calC$  of minimal Van Kampen diagrams in $\tilHi$ for the words in  $\calV_i'$. Since every triangular relation of $\tilHi$ is allowed, this amounts to find minimal triangulations of a polygon.

        \item  Make the list $\calD$ of all Van Kampen diagrams with at most $240K$ cells in $\calR$ or in a diagram in the list $\calC$. Since the set of permited cells, and their numbers are finite and known,   this list is finite, and can be found by enumeration.

        \item  Make $\calW$ the list of words labeling the boundary of a diagram in $\calD$.

        \item  For each element $w$ of the (finite) list $\calW$, find the subset of $\calD$ of diagrams whose boundaries are labeled by $w$, and find a smallest element (for the area) in this subset, and finally compute $A(w)$ the area of this smallest element.

        \item  For all $w \in \calW$, such that  $A(w) \leq 240 K$,
          check whether $A(w)> (\sqrt{K}/600) |w\absx$.  If this is true for
          some $w$, then $K:=K+1$ and go back to the first instruction, and if
          this is false for every $w$, then stop the algorithm, and return the integer $K$. 
 
        \end{enumerate}

We  assume that either the presentation is exact, or a solution to the word problem of $\Gamma$ was given as input.

        Either the algorithm stops at step $1$ for some $K$, or it stops at step $7$ for some $K$, or it never stops.   

        If the  algorithm  stops at step $1$, this is because it has found a relation of $\Gamma$ with letters in $\tilHi$ that is not a relation of $\tilHi$, and therefore, the presentation was indeed non exact.

        If the  algorithm does not detect that the presentation is non-exact, it does not stop at step $1$ 
        (for instance because the presentation is indeed exact).   
        Then, for all $K$, the assumption of  Lemma \ref{lem;reduction_finale} is still satisfied.    
        In this case, by Lemma \ref{lem;find_minimal}, 
         for all $w\in \calW$, $A(w)=Area(w)$.

        If the presentation satisfies a linear isoperimetric inequality with factor $K_0$, then the algorithm will stop (at step 7) for some $K\leq (600 K_0)^2$. 

	Conversely, if the algorithm stops at step 7 for $K$, then every relation $w''$ with a minimal 
        Van Kampen diagram in $\calD$  
        satisfies  $Area(w'')\leq \sqrt{K}/600 |w''\absx$, and therefore by 
        Lemma \ref{lem;reduction_finale},    for all  relation $w$, $Area(w) \leq K |w\absx$. 
        The presentation satisfies a linear isoperimetric inequality of factor $K$.

        This ensures that if no isoperimetric relation is satisfied, and if the algorithm 
        does not stop at step $1$, it will never stop.

Assume, finally, that the presentation provided is non exact, but a solution to the word problem of $\Gamma$ was given in input.  Assume moreover that some quotient $\tilHi \tto H_i$ has infinite kernel. There is a sequence  $w_k$  of different letters in $\tilHi$ in the kernel. They produce an infinite sequence of relations of $\Gamma$ of relative length $1$, and   by Lemma \ref{lem;find_minimal}, their  areas cannot remain bounded.  

        Thus the presentation satisfies no isoperimetric inequality, and as
        proved before, the algorithm cannot stop at step $7$. Therefore it
        will stop at step $1$ for some  $K$ smaller than the length of the
        element $w_1$ for the word metric induced by $\calX_\calR \cap
        \tilHi$.  \end{proof}

    Let us recall that Farb's coned-off cayley graph $\widehat{Cay}\Gamma$ of a relatively hyperbolic
    group $\Gamma$ 
    is defined by the following construction. Take the Cayley graph
    on the given generating set, and add  one 
    vertex  
    $v_{[\gamma H_k]}$  
    for each left coset $[\gamma H_k]$ of each parabolic subgroup
    $H_k$, and one edge $(v_{[\gamma H_k]}, \gamma h )$, for each element $
    h\in H_k$.  If $\Gamma$ is relatively hyperbolic, 
    this is a hyperbolic graph, acted upon by
    $\Gamma$, and although it is non locally finite in general (when some
    $H_k$ is infinite), there are only finitely many orbits of 
    simple loops of given
    length (see \cite{F}, \cite{Brel} for details).

\begin{coro}\label{coro;delta}
  There is an algorithm that, given an exact finite relative presentation of a
  relatively hyperbolic group, provides an hyperbolicity constant for
  $\widehat{Cay}\Gamma$, and the consecutive lists of orbits of simple loops
  in  $\widehat{Cay}\Gamma$ of given length.
\end{coro}

 \begin{proof} An hyperbolicity constant of the Cayley graph 
	with set of generators
  $\calX$, can be computed from an 
      explicit linear
      Dehn function (see \cite[2.9 p.~419]{BH}, or \cite[p.~61]{Rebb} 
	(in which constants are
      given at the end of the argument page 69-70).  
      Therefore, one can compute such a constant 
      for the conned-off  graph on $X$, since the identity on the vertices of
      $Cay \Gamma$ is 2-bi-Lipshitz.

        The list of simple loops of given length in  $\widehat{Cay}\Gamma$ is
        easily seen to
        correspond to the list of relations of given relative length, whose
        Van Kampen diagram do not have a cluster 
        with two edges on the boundary. To make the list of those of given
        length, it is enough to make the list of those of given area, but
        according to Lemma \ref{lem;find_minimal}, one only needs 
        to use cells with
        perimeter labeled by letters of absolute length bounded by an explicit
        constant. 
        This is enough to enumerate effectively all possible candidates, and
        a comparison with the list of shorter simple loops indicates which are
        actually simple.  \end{proof}

\subsection{Abelian parabolics: an algorithm to find them.}

In general, given a finite subset of a group, it is hard to tell what kind of subgroup it generates. But if the elements commute pairwise, they generate an abelian subgroup. 
This motivates the sharper study of 
      relatively hyperbolic groups,
      with  abelian parabolic subgroups. We start with an arbitrary
      presentation of such a group, and we want to find where are the
      parabolic subgroups, and an exact relative presentation.

The word problem is solvable for these groups, but some solutions require the preliminary knowledge of the hyperbolicity constant, and of the ranks of
      the abelian subgroups, which we don't have {\it a priori}, since we
      don't have an exact relative presentation yet. But fortunately, we can proceed without this assumption, as pointed out to us by M.~Bridson.

      Recall that, by a theorem of  D. Rebbechi \cite{Rebb}, any relatively hyperbolic group with  abelian subgroups 
      is bi-automatic. Thus one can use first the following, derived from the algorithm of Todd-Coxeter, and the quadratic time solution to the word problem for automatic groups:

      \begin{prop}\label{lem;auto}(\cite[Theorems 5.2.4 and 2.3.10]{E})
 
        There is an algorithm whose input is a finite presentation of a 
        group, that stops if and only if the group is automatic, and that outputs a solution to the word problem.

      \end{prop}

    \begin{proof}[of Theorem  \ref{theo;intro2bis}.]
      Let $\Gamma$ be a  group with a finite presentation 
      $\langle X|R \rangle$, and assume that
      $\Gamma$ is hyperbolic relative to abelian subgroups.

      Enumerate the triples of integers  $(n, k, s)$, 
      and for each of them, the families of
      $n$ subsets $\pi_1,\dots, \pi_n$ of at most $k$ elements of $\Gamma$,
      such that    
      for all $i$, all the elements in $\pi_i$
      commute pairwise in $\Gamma$ 
      (for that we use a solution of the word problem in $\Gamma$, 
        which we have, by Proposition
      \ref{lem;auto}). 

        Given such a family, each $\pi_i$ generates an 
        abelian subgroup $H_i$ of $\Gamma$.  
        For each element of $\pi_i$, check with a solution to the word problem
        whether it is of order $s$ or less in $\Gamma$.

        One can then compute a (not necessarily exact)
      relative presentation of $\Gamma$, by considering $\tilHi$ 
      the 
      abelian group of basis $\pi_i$, whose every generator has either
      infinite order, or order at most $s$ as found by the latter verification.

      By structure of finitely generated abelian groups, one of the relative
      presentations above (for some correct $(n, k, s)$ and $\pi_i $)
      corresponds to an exact relative presentation of $\Gamma$.

      In order to finish the proof, it suffices to recognize when
      such a presentation is exact, and corresponds to the relative structure,  
      that is, when the sets $\pi_i$ are
      actually basis in $\Gamma$ of maximal abelian subgroups of
      $\Gamma$, with correct orders.

        For each $i$, if the  map $\tilHi \tto H_i$ has non trivial but finite
        kernel, then it lies in the factor of $\tilHi$ generated by the
        elements of $\tilHi$ of finite order. This implies that one of them is
        in $\Gamma$ of order less than $s$, but we assumed that we checked 
        that it
        does not happen.  Therefore the  map $\tilHi \tto H_i$  
        is either injective, or has infinite kernel.
        The algorithm of Proposition \ref{prop;algo}, with, as input, the proposed presentation, and a solution to the word problem for $\Gamma$, will then stop if the presentation is non-exact, indicating this fact. On the other hand, if the presentation is exact, it will stop if and only if the presentation satisfies a linear isoperimetric inequality. In such case, the presentation corresponds to the structure of a relatively hyperbolic group. \end{proof}

% FIN DE LA REVISION

\bigskip

   {\sc 
     Institut de Math\'ematiques de Toulouse
	Laboratoire E. Picard, 
     Univ. P. Sabatier, 
     F-31062 Toulouse, France.
   }

   {\it 
     E-mail:  dahmani@math.univ-toulouse.fr
   }

\end{document}